\documentclass[12pt]{article} %***
\usepackage[sectionbib]{natbib}
\usepackage{array,epsfig,fancyheadings,rotating}
\usepackage[]{hyperref}  
\usepackage{tikz,tikz-3dplot}
\usepackage{sectsty, secdot}
\sectionfont{\fontsize{12}{14pt plus.8pt minus .6pt}\selectfont}
\renewcommand{\theequation}{\arabic{equation}}
\subsectionfont{\fontsize{12}{14pt plus.8pt minus .6pt}\selectfont}
\usepackage{amsmath}
\usepackage{amssymb}
\usepackage{amsfonts}
\usepackage{multirow}
\usepackage{amsthm}
\usepackage{version}
\usepackage{graphics}

\newtheorem{theorem}{Theorem}
\newtheorem{lemma}{Lemma}

\newtheorem{proposition}{Proposition}
\theoremstyle{definition}

\newcommand{\cqfd}{\hfill $\square$}

\newcommand{\R}{\mathbb R}

\newcommand{\n}{^{(n)}}

\newcommand{\sirc}{{\scriptscriptstyle
\circ }}
\newcommand{\vech}{{\rm ve}\!^{^{\sirc}}\!\!{\rm ch}\,}
\newcommand{\Xb}{\mathbf{X}}

\newcommand{\Sb}{\mathbf{S}}

\newcommand{\Vb}{\mathbf{V}}

\newcommand{\Zb}{\mathbf{Z}}

\newcommand{\xb}{\ensuremath{\mathbf{x}}}

\newcommand{\Ab}{\ensuremath{\mathbf{A}}}

\newcommand{\Db}{\ensuremath{\mathbf{D}}}

\newcommand{\Yb}{\ensuremath{\mathbf{Y}}}

\newcommand{\thetab}{{\pmb \theta}}

\newcommand{\mub}{{\pmb \mu}}
\newcommand{\Umb}{{\pmb \Upsilon}}
\newcommand{\Thetab}{{\pmb \Theta}}
\newcommand{\Lamb}{{\pmb \Lambda}}
\newcommand{\Omegab}{{\pmb \Omega}}
\newcommand{\Sigb}{{\pmb \Sigma}}
\newcommand{\deltab}{{\pmb \delta}}
\newcommand{\Deltab}{{\pmb \Delta}}
\newcommand{\taub}{{\pmb \tau}}
\newcommand{\nub}{{\pmb \nu}}

\newcommand{\Gamb}{{\pmb \Gamma}}

\newcommand{\sigbu}{\underline{\pmb \sigma}}

\newcommand{\I}{{\scriptscriptstyle I}}
\newcommand{\II}{{\scriptscriptstyle I\!I}}

\newcommand{\pr}{^{\prime}}

\newcommand{\ny}{n\rightarrow\infty}

\renewcommand{\theequation}{\thesection.\arabic{equation}}

\begin{document}

%%%%%%%%%%%%%%%%%%%%%%%%%%%%%%%%%%%%%%%%%%%%%%%%%%%%%%%%%%%%%%%%%%%%%%%%%%%%%%%%%%%%%%%%%%%%%%%%%%%%%%%%%%%%%%%%%%%%%%%%%%%%
%%%%%%%%%%%%%%%%%%%%%%%%%%%%%%%%%%%%%%%%%%%%%%%%%%%%%%%%%%%%%%%%%%%%%%%%%%%%%%%%%%%%%%%%%%%%%%%%%%%%%%%%%%%%%%%%%%%%%%%%%%%%

\renewcommand{\baselinestretch}{1.3}

\markright{ \hbox{\footnotesize\rm %Statistica Sinica
%{\footnotesize\bf 24} (201?), 000-000
}\hfill\\[-13pt]
\hbox{\footnotesize\rm
%\href{http://dx.doi.org/10.5705/ss.20??.???}{doi:http://dx.doi.org/10.5705/ss.20??.???}
}\hfill }

\markboth{\hfill{\footnotesize\rm Davy Paindaveine, Jos\'ea  Rasoafaraniaina and Thomas Verdebout} \hfill}
{\hfill {\footnotesize\rm Preliminary test estimation in ULAN models} \hfill}

\renewcommand{\thefootnote}{}
$\ $\par

%%%%%%%%%%%%%%%%%%%%%%%%%%%%%%%%%%%%%%%%%%%%%%%%%%%%%%%%%%%%%%%%%%%%%%%%%%%%%%%%%%%%%%%%%%%%%%%%%%%%%%%%%%%%%%%%%%%%%%%%%%%%

\fontsize{12}{14pt plus.8pt minus .6pt}\selectfont \vspace{0.8pc}
\centerline{\large\bf Preliminary Test Estimation in ULAN Models}
%\vspace{2pt}
%\centerline{\large\bf Locally Asymptotically Normal models}
\vspace{.4cm}
\centerline{Davy Paindaveine$^{a,b}$, Jos\'ea  Rasoafaraniaina$^a$ and Thomas Verdebout$^a$}
\vspace{.4cm}
\centerline{\it $^a$Universit\'e libre de Bruxelles (ULB)}
\centerline{\it $^b$Toulouse School of Economics (TSE)}
%\vspace{.2cm}
%\centerline{\it $^\dagger$Universit\'e libre de Bruxelles (ULB)}
\vspace{.55cm}
\fontsize{9}{11.5pt plus.8pt minus .6pt}\selectfont

%%%%%%%%%%%%%%%%%%%%%%%%%%%%%%%%%%%%%%%%%%%%%%%%%%%%%%%%%%%%%%%%%%%%%%%%%%%%%%%%%%%%%%%%%%%%%%%%%%%%%%%%%%%%%%%%%%%%%%%%%%%%

\begin{quotation}
\noindent {\it Abstract:} Preliminary test estimation, which is a natural procedure when it is suspected a priori that the parameter to be estimated might take value in a submodel of the model at hand, is a classical topic in estimation theory. In the present paper, we establish general results on the asymptotic behavior of preliminary test estimators. More precisely, we show that, in uniformly locally asymptotically normal (ULAN) models, a general asymptotic theory can be derived for preliminary test estimators based on estimators admitting generic Bahadur-type representations. This allows for a detailed comparison between  classical estimators and preliminary test estimators in ULAN models. Our results, that, in standard linear regression models, are shown to reduce to some classical results, are also illustrated in more modern and involved setups, such as the multisample one where $m$ covariance matrices~${\pmb\Sigma}_1, \ldots, {\pmb\Sigma}_m$ are to be estimated when it is suspected that these matrices might be equal, might be proportional, or might share a common ``scale". Simulation results confirm our theoretical findings.

\par

\vspace{9pt}
\noindent {\it Key words and phrases:}
LAN models, Le Cam's asymptotic theory, Multisample covariance matrix estimation, Preliminary test estimation.
\par
\end{quotation}\par

\renewcommand{\theequation}{\thesection.\arabic{equation}}

\def\thefigure{\arabic{figure}}
\def\thetable{\arabic{table}}

\fontsize{12}{14pt plus.8pt minus .6pt}\selectfont

%\newpage

\setcounter{section}{1}
\setcounter{equation}{0} %-1
\noindent {\bf 1. Introduction}
\vspace{4mm}

\noindent Preliminary test estimation is a widely studied topic in Statistics and Econometrics, that can be traced back to the seminal paper by \cite{Ban44}. Preliminary test estimators are typically useful when one has to perform statistical inference with some ``uncertain prior information".  More formally, assume that one is interested in estimating a parameter~$\thetab$ that belongs to some parameter space $\Thetab \subset \R^p$, with the  ``uncertain prior information" that~$\thetab$ belongs to a given subset $\Thetab_0$ of $\Thetab$ (throughout, we assume that $\Thetab$ is an open subset of $\R^p$). Then, roughly speaking, the statistician may hesitate between (i) an \emph{unconstrained} estimator~$\hat{\thetab}_{\rm U}$ with values in~$\Thetab$ or 
\linebreak
(ii) a \emph{constrained} estimator~$\hat{\thetab}_{\rm C}$ with values in~$\Thetab_0$ only. The idea underpinning preliminary test estimation is relatively simple: if a suitable test~$\phi_n$ for~${\cal H}_0: \thetab \in \Thetab_0$ against~${\cal H}_1: \thetab \notin \Thetab_0$ did not reject the null hypothesis, then~$\hat{\thetab}_{\rm C}$ should be used; on the contrary, if~$\phi_n$ provided evidence against~${\mathcal H}_0$, then the unconstrained estimator~$\hat{\thetab}_{\rm U}$ should be favoured. A preliminary test estimator based on the estimators $\hat{\thetab}_{\rm U}$ and $\hat{\thetab}_{\rm C}$ and on the test $\phi_n$ is therefore
\begin{equation} 
\label{preliintro}
\hat{\thetab}_{\rm PTE}:= 
{\mathbb I}[\phi_n=1] \hat{\thetab}_{\rm U}
+
{\mathbb I}[\phi_n=0] \hat{\thetab}_{\rm C}
,
\end{equation}
where $ {\mathbb I}[A]$ stands for the indicator function associated with~$A$ and where~$\phi_n=1$ (resp., $\phi_n=0$) indicates rejection (resp., non-rejection) of~${\mathcal H}_0$ by~$\phi_n$.  
%The idea underpinning $\hat{\thetab}_{\rm PTE}$ in \eqref{preliintro} is simply to take $\hat{\thetab}_{\rm U}$ is the test $\phi$ rejects ${\cal H}_0$ and $\hat{\thetab}_{\rm C}$ if $\phi$  does not reject ${\cal H}_0$.

Since \cite{Ban44}, preliminary test estimation has been an active research topic. \cite{SeSa79}, \cite{SeSa87}, \cite{WA06} and \cite{KS14} considered preliminary test estimation in regression models. \cite{GLG92} tackled the problem of selecting the size of the test~$\phi_n$ when conducting preliminary test estimation in a misspecified regression model. \cite{OT80} considered estimation of regression coefficients after a preliminary test for homoscedasticity. Preliminary test estimation in elliptical models has been considered in \cite{AA14} and by \cite{PRV17} in a principal component analysis context. Preliminary test estimation has also been widely considered in time series models; see, e.g., \cite{AB00}, \cite{Tani11}, and the references therein. For a general overview of the topic, we refer to \cite{GG93} and \cite{SA06}.

%In most of the contributions above, preliminary test estimation is considered in regular models. 
Despite the many works on the topic, there does not seem to exist a general theory describing the asymptotic behavior of preliminary test estimators. The main objective of the present paper is therefore to derive such a general theory and to do so  in a broad class of models (that will include in particular all models mentioned above). Assuming that the underlying model is regular in the sense that it is uniformly locally asymptotically normal (ULAN), we will derive the asymptotic behavior of a general preliminary test estimator; more precisely, we will consider preliminary test estimators based on estimators~$\hat{\thetab}_{\rm U}$ and~$\hat{\thetab}_{\rm C}$ that admit Bahadur-type representations. Our asymptotic results do cover many of the existing results in the literature but also allow us to consider more modern and involved models.  

As expected, the asymptotic behavior of preliminary test estimators will depend on the true value of the parameter~$\thetab$. We first show that when this true value is fixed outside~$\Thetab_0$, then,  provided that the test~$\phi_n$ is consistent, a preliminary test estimator is asymptotically equivalent in probability to the unconstrained estimator $\hat{\thetab}_{\rm U}$. Second, we show that when the true value of~$\thetab$ asymptotically belongs to contiguous regions of~$\Thetab_0$ (in a sense that is related to the asymptotic concept of contiguity, as we will make precise below), a preliminary test estimator exhibits an asymptotic behavior achieving a nice compromise between~$\hat{\thetab}_{\rm U}$ and~$\hat{\thetab}_{\rm C}$. 

The paper is organized as follows. In Section~2, we describe the assumptions that will be considered in the sequel. In Section 3, we state our asymptotic results and derive explicit forms for the asymptotic mean square error of preliminary test estimators based on asymptotically efficient estimators. In Section 4, we illustrate these general results in two particular setups. First, we show that, in a simple linear regression context, our results allow us to recover the classical results from \cite{SA06}. Then, we consider preliminary test estimation of $m$ covariance matrices in a multisample Gaussian setup. Preliminary test estimators associated with the constraints of \emph{covariance homogeneity}, \emph{shape homogeneity} and \emph{scale homogeneity} are studied. Monte Carlo simulations confirm our theoretical results. Finally, an appendix collects the proofs.
\vspace{4mm}

%%%%%%%%%%%%%%%%

\newpage
\setcounter{section}{2}
\setcounter{equation}{0} %-1
\noindent {\bf 2. ULAN models and Preliminary Test Estimators}
\vspace{4mm}

\noindent As mentioned in the introduction, our objective is to derive the asymptotic behavior of preliminary test estimators (PTEs) in a very general context. 
We will throughout assume that the underlying parametric model $\{ {\rm P}_\thetab\n: \thetab \in \Thetab \subset \R^p\}$ under investigation is \emph{uniformly locally and asymptotically normal (ULAN)} in the following sense (throughout, all convergences are as~$n\to\infty$). 
\vspace{2mm}

\noindent {\bf Assumption~A}. There exists a sequence $(\nub_n)$ of full-rank non-random $p\times p$ matrices that is~$o(1)$ and a sequence~$(\thetab_n)$ in~$\Thetab$ with $\nub_n^{-1}(\thetab_n- \thetab)=O(1)$ for some $\thetab \in \rm \Thetab$ such that for every sequence~$(\taub_n)$ that is~$O(1)$ and satisfies~$\thetab_n+ \nub_n \taub_n \in \Thetab$ for any~$n$,
\begin{equation}
\Lambda\n:=\log \frac{{\rm dP}^{(n)}_{\thetab_n+\nub_n \taub_n}}{{\rm d P}^{(n)}_{\thetab_n}} = \taub_n\pr \pmb{\Delta}^{(n)}_{\thetab_n}-\frac{1}{2} \taub_n\pr\Gamb_\thetab \taub_n +o_{{\rm P}}(1)
\end{equation}
under ${\rm P}^{(n)}_{\thetab}$, where the random $p$-vector~$\pmb{\Delta}^{(n)}_{\thetab}$, still under~${\rm P}\n_\thetab$, is asymptotically normal with mean vector~${\bf 0}$ and covariance matrix $ \Gamb_{\thetab}$.
% is such that there exists an \mbox{i.i.d.} sequence ${\bf A}_{1, \thetab} , \ldots, {\bf A}_{n, \thetab} $ with ${\rm E}[{\bf A}_{1, \thetab} ]={\bf 0}$ and ${\rm Var}[{\bf A}_{1, \thetab} ]= \Gamb_{\thetab} <+ \infty$ such that $\Deltab_{\thetab}\n = \frac{1}{\sqrt{n}} \sum_{i=1}^n {\bf A}_{i, \thetab}  + o_{\rm P}(1)$ still under~${\rm P}\n_\thetab$ .
%asymptotically normal with mean zero and covariance matrix $\Gamb_\thetab <+\infty$ still under ${\rm P}^{(n)}_{\thetab}$ as $n \to \infty$.
\vspace{2mm}

An extensive list of models do satisfy Assumption~A. This list includes hidden Markov models (\citealp{BR96}), quantum mechanics models (\citealp{KG09}, \citealp{GK15}), time series models (\citealp{DKW97}, \citealp{HTSC99}, \citealp{FZ13}), elliptical models (\citealp{HP06}, \citealp{HPV10}), multisample elliptical models (\citealp{HP08}, \citealp{HPV13}, \citealp{HPV14}), models for directional data (\citealp{LSTV13}, \citealp{GPV19}), and many more.

As explained in the introduction, the construction of a PTE involves an unconstrained estimator~${\hat \thetab}_{\rm U}$ taking values in~$\Thetab$, a constrained estimator~${\hat \thetab}_{\rm C}$  taking values in~$\Thetab_0$, and a test~$\phi_n$ for~$\mathcal{H}_0:\thetab\in\Thetab_0$ against~$\mathcal{H}_1:\thetab \notin\Thetab_0$. Throughout, we will assume that~$\Thetab_0$ is a linear subspace of~$\R^p$ of the form 
$$
\Thetab_0 = (\thetab_0 + {\mathcal M}(\Umb)) \cap \Thetab,
$$
where~$\thetab_0$ is a fixed $p$-vector and ${\mathcal M}(\Umb)$ denotes the vector subspace of~$\R^p$ that is spanned by the columns of the $p \times r$ full-rank matrix $\Umb$ ($r< p$). We will restrict to the case~$\thetab_0={\bf 0}$, which is without loss of generality (a reparametrization of the model always allows us to reduce to this case). Throughout, we will consider PTEs of the form
$$
\hat{\thetab}_{\rm PTE}
=
 {\mathbb I}[\phi_n=1] \hat{\thetab}_{\rm U}
 +
 {\mathbb I}[\phi_n=0] \hat{\thetab}_{\rm C}
$$ 
that involve estimators~$\hat{\thetab}_{\rm U}$, $\hat{\thetab}_{\rm C}$ and a test~$\phi_n$ satisfying the following assumption.
\vspace{2mm}

\noindent {\bf Assumption~B}. With~$\nub_n$ and $\Deltab_\thetab\n$ as in Assumption~A, there exists a random $p$-vector ${\bf S}_\thetab\n$ for which
$({\bf S}_\thetab^{(n)\prime}, \Deltab_\thetab^{(n)\prime})\pr$ is asymptotically normal with mean vector ${\bf 0}$ and covariance matrix 
$$
\Bigg(
\begin{array}{cc} 
\Sigb_\thetab & \Omegab_\thetab \\[-2mm]  
\Omegab_\thetab & \Gamb_\thetab \end{array}
\Bigg)
$$
 under~${\rm P}\n_\thetab$ and such that, for some $p \times p$ matrix~${\bf A}_\thetab$ and~$r \times p$ matrix~${\bf B}_\thetab$, 
\begin{itemize}
\item[(i)] 
 $\nub_n^{-1} (\hat \thetab_{\rm U}- \thetab)= {{\bf A}}_\thetab {\bf S}_\thetab\n + o_{\rm P}(1)$ under~${\rm P}\n_\thetab$, $\thetab \in \Thetab$,
\item[(ii)] $\nub_n^{-1} (\hat \thetab_{\rm C}- \thetab)= \Umb { \bf B}_\thetab {\bf S}_\thetab\n+ o_{\rm P}(1)$
 under~${\rm P}\n_\thetab$, $\thetab \in \Thetab_0$, and
\item[(iii)] $\phi_n$ rejects~$\mathcal{H}_0:\thetab\in\Thetab_0$ at asymptotic level $\alpha$ when 
$ {Q}^{(n)} : = \| {\bf D}\n \|^2 > \chi^2_{p-r,1- \alpha},$
where $\chi^2_{\ell,\beta}$ denotes the upper $\beta$-quantile of the~$\chi^2_\ell$ distribution and where ${\bf D}\n$ is such that ${\bf D}\n= {\bf C}_\thetab {\bf S}_\thetab\n + o_{\rm P}(1)$  under~${\rm P}\n_\thetab$, $\thetab \in \Thetab_0$, for some $p \times p$ matrix
${\bf C}_\thetab$ satisfying $\Sigb_\thetab{\bf C}_\thetab\pr {\bf C}_\thetab \Sigb_\thetab{\bf C}_\thetab\pr  {\bf C}_\thetab \Sigb_\thetab= \Sigb_\thetab{\bf C}_\thetab\pr {\bf C}_\thetab \Sigb_\thetab$ and ${\rm tr}[{\bf C}_\thetab\pr {\bf C}_\thetab \Sigb_\thetab]=p-r$. Furthermore, ${\rm P}[{Q}^{(n)}>\chi^2_{p-r,1- \alpha}]$ converges to one under~${\rm P}\n_\thetab$, $\thetab \notin \Thetab_0$.
\end{itemize}
\vspace{2mm}

While being quite complex, Assumption~B is extremely mild and, provided that the underlying model is ULAN as in Assumption~A, merely only imposes that the unconstrained estimator~$\hat{\thetab}_{\rm U}$ admits a Bahadur-type representation. To show this, restrict to the usual contiguity rate~$\nub_n= n^{-1/2}{\bf I}_p$ (extension to a general~$\nub_n$ is direct) and assume that, under~${\rm P}\n_\thetab$, $\thetab\in \Thetab$,   
\begin{equation} 
\label{Bahadurrep}
\sqrt{n}(\hat \thetab_{\rm U}- \thetab)= \frac{1}{\sqrt{n}} \sum_{i=1}^n {\bf T}_i\n +o_{\rm P}(1),
\end{equation}
 where the random $p$-vectors~${\bf T}_i\n={\bf T}_i\n(\thetab)$, $i=1,\ldots,n$ are mutually independent and share a common distribution that has mean zero and has finite second-order moments (this ensures that Assumption~B(i) holds, with~${\bf B}_\thetab:=\mathbf{I}_p$ and~${\bf S}_\thetab^{(n)}:=n^{-1/2} \sum_{i=1}^n {\bf T}_i\n$, say). Under very mild assumptions (needed to check the Levy-Lindeberg condition), the CLT for triangular arrays will then ensure that~$({\bf S}_\thetab^{(n)\prime}, \Deltab_\thetab^{(n)\prime})\pr$ is asymptotically normal under~${\rm P}\n_\thetab$, as required in Assumption~B.   
Letting~${\bf P}_\Umb:=\Umb (\Umb\pr \Umb)^{-1} \Umb\pr$ be the matrix of the projection onto the constraint~$\Thetab_0={\cal M}(\Umb)\cap\Thetab$, the constrained estimator~${\hat \thetab}_{\rm C} := {\bf P}_\Umb {\hat \thetab}_{\rm U}$ readily satisfies 
$$
\sqrt{n}(\hat \thetab_{\rm C}- \thetab)= {\bf P}_\Umb \sqrt{n}(\hat \thetab_{\rm U}- \thetab)=  \frac{1}{\sqrt{n}} \sum_{i=1}^n {\bf P}_\Umb{\bf T}_i\n +o_{\rm P}(1)
$$
under~${\rm P}_\thetab\n$, $\thetab \in \Thetab_0$, so that Assumption~B(ii) is fulfilled, too (with~${\bf B}_\thetab:=(\Umb\pr \Umb)^{-1} \Umb\pr$). 
Finally, note that Assumption~B(iii) will be satisfied by Wald tests for~${\cal H}_0: \thetab \in \Thetab_0$ against ${\cal H}_1: \thetab \notin \Thetab_0$  constructed in the usual way from~\eqref{Bahadurrep}. Wrapping up, the only key point in Assumption~B is its part~(i), which itself holds as soon as the unconstrained estimator~$\hat \thetab_{\rm U}$, like, e.g., most M-, R-, and S-estimators, admits a Bahadur-type representation. 

Now, in the ULAN framework of Assumption~A, it should be noted that an asymptotically efficient (unconstrained) estimator~$\hat{\thetab}_{\rm U}$, that is, an estimator that, under~${\rm P}\n_\thetab$, $\thetab \in \Thetab$, is such that 
\begin{equation} 
\label{optiU}
\sqrt{n}(\hat \thetab_{\rm U}- \thetab)= \Gamb_\thetab^{-1} \Deltab_\thetab\n+ o_{\rm P}(1) 
\end{equation}
(see, e.g., Chapter~3 of \citealp{TK00}) also satisfies Assumption~B(i), with~${\bf A}_\thetab=\Gamb_\thetab^{-1}$ and $\Sb_\thetab\n= \Deltab_\thetab\n$ (which provides~$\Sigb_{\thetab}=\Gamb_{\thetab}$).
An asymptotically efficient constrained estimator~$\hat{\thetab}_{\rm C}$, that is such that 
\begin{equation} \label{optiC}
\sqrt{n}(\hat \thetab_{\rm C}- \thetab)= \Umb (\Umb\pr \Gamb_\thetab \Umb)^{-1} \Umb\pr \Deltab_\thetab\n + o_{\rm P}(1)
\end{equation}
 under~${\rm P}\n_\thetab$, $\thetab \in \Thetab_0$,
 satisfies Assumption~B(ii), with ${\bf B}_\thetab=(\Umb\pr \Gamb_\thetab \Umb)^{-1} \Umb\pr$ and $\Sb_\thetab\n= \Deltab_\thetab\n$.  For testing ${\cal H}_0: \thetab \in \Thetab_0$ against ${\cal H}_1: \thetab \notin \Thetab_0$, the locally asymptotically most stringent test rejects~${\cal H}_0$ at asymptotic level~$\alpha$ when
\begin{equation} \label{optitest}
Q\n
= 
\big\| 
({\bf I}_p- \Gamb_{\hat \thetab_{\rm C}}^{1/2}\Umb(\Umb\pr \Gamb_{\hat \thetab_{\rm C}} \Umb)^{-1}\Umb\pr\Gamb_{\hat \thetab_{\rm C}}^{1/2}) \Gamb_{\hat \thetab_{\rm C}}^{-1/2} \Deltab_{\hat \thetab_{\rm C}}\n 
\big\|^2
> 
\chi^2_{p-r,1- \alpha};
\end{equation}
see, e.g., Chapter 5 of \cite{LV17}. Under Assumption~A, it is easy to check that, under~${\rm P}\n_\thetab$, $\thetab \in \Thetab_0$,
\begin{eqnarray*}
	\lefteqn{
({\bf I}_p- \Gamb_{\hat \thetab_{\rm C}}^{1/2}\Umb(\Umb\pr \Gamb_{\hat \thetab_{\rm C}} \Umb)^{-1}\Umb\pr\Gamb_{\hat \thetab_{\rm C}}^{1/2}) \Gamb_{\hat \thetab_{\rm C}}^{-1/2} \Deltab_{\hat \thetab_{\rm C}}\n
}
\\[2mm]
 & & 
 \hspace{13mm} 
 =
 ({\bf I}_p- \Gamb_{\thetab}^{1/2}\Umb(\Umb\pr \Gamb_{\thetab} \Umb)^{-1}\Umb\pr\Gamb_{\thetab}^{1/2}) \Gamb_{\thetab}^{-1/2} \Deltab_{\thetab}\n+ o_{\rm P}(1)
,
\end{eqnarray*}
so that Assumption~B(iii) then holds, still with~$\Sb_\thetab\n= \Deltab_\thetab\n$, $\Sigb_{\thetab}=\Gamb_{\thetab}$, and with 
$
{\bf C}_\thetab= ({\bf I}_p- \Gamb_{\thetab}^{1/2}\Umb(\Umb\pr \Gamb_{\thetab} \Umb)^{-1}\Umb\pr \Gamb_{\thetab}^{1/2}) \Gamb_{\thetab}^{-1/2}=({\bf I}_p- \Gamb_{\thetab}^{1/2}\Umb {\bf B}_\thetab \Gamb_{\thetab}^{1/2}) \Gamb_{\thetab}^{-1/2}
$ 
(one can indeed check that~${\bf C}_\thetab\pr {\bf C}_\thetab \Gamb_\thetab {\bf C}_\thetab\pr {\bf C}_\thetab= {\bf C}_\thetab\pr {\bf C}_\thetab$ and that
${\rm tr}[{\bf C}_\thetab\pr {\bf C}_\thetab \Gamb_\thetab]={\rm tr}[{\bf I}_p]- {\rm tr}[(\Umb\pr \Gamb_{\thetab} \Umb)^{-1}(\Umb\pr \Gamb_{\thetab} \Umb)]=p-r$).
%%\begin{eqnarray}
%%{Q}^{(n)} &=&  \| \Sigb_{{\hat \thetab}_{\rm U}}^{1/2} (\nub\n)^{-1}(\hat \thetab_{\rm U} - \hat \thetab_{\rm C})\|^2 \nonumber \\
%%&=&  \| \Sigb_{\thetab}^{1/2}(\Sigb_\thetab^{-1}-\Umb(\Umb\pr \Sigb_\thetab \Umb)^{-1}\Umb\pr)  \Deltab_\thetab\n  \|^2+ o_{\rm P}(1) \nonumber \\
%%&=& \| ({\bf I}_p- \Sigb_{\thetab}^{1/2}\Umb(\Umb\pr \Sigb_\thetab \Umb)^{-1}\Umb\pr\Sigb_{\thetab}^{1/2}) \Sigb_{\thetab}^{-1/2} \Deltab_\thetab\n  \|^2+ o_{\rm P}(1)
%%\end{eqnarray}
%% under~${\rm P}\n_\thetab$ with $\thetab \in \Thetab_0$. In the sequel we use the notations  ${\bf P}_\Umb={\bf I}_p - \Sigb_\thetab^{1/2} {\bf A}_\Umb \Sigb_\thetab^{1/2},$ 
%%and ${\bf A}_\Umb=\Umb (\Umb\pr \Sigb_\thetab \Umb)^{-1}\Umb\pr$ for the sake of readibility. Using these notations, we have that $ {Q}^{(n)}=\| {\bf D}\n \|^2$ where
%%\begin{eqnarray} \label{Dasymp}
%%{\bf D}\n &= &\Sigb_{{\hat \thetab}_{\rm U}}^{1/2} (\nub\n)^{-1}(\hat \thetab_{\rm U} - \hat \thetab_{\rm C}) \nonumber \\
%%&= & ({\bf I}_p - \Sigb_\thetab^{1/2} {\bf A}_\Umb \Sigb_\thetab^{1/2}) \Sigb_\thetab^{-1/2}{\bf S}_\thetab\n+o_{\rm P}(1) \nonumber \\
%%& = & {\bf P}_\Umb  \Sigb_\thetab^{-1/2} {\bf S}_\thetab\n+o_{\rm P}(1)  
%%\end{eqnarray}
%%still  under~${\rm P}\n_\thetab$ with $\thetab \in \Thetab_0$.
To summarize, Assumptions~A and B cover many existing models and estimators. In the next section, our objective is to derive asymptotic results for PTEs in the general framework covered by these assumptions.
\vspace{4mm}

%%%%%%%%%%%%%%%%%%%%%%%%%%%%%%%%%%%%%%%%%%%%

\setcounter{section}{3}
\noindent {\bf 3. Asymptotic results}
\vspace{4mm}

\noindent In this section, we derive the asymptotic behavior of a PTE of the form 
\begin{equation} \label{PTEinvest}
{\hat \thetab}_{\rm PTE}:= {\mathbb I}[\phi_n=1] {\hat \thetab}_{\rm U}+ {\mathbb I}[\phi_n=0] {\hat \thetab}_{\rm C},
\end{equation}
where the estimators ${\hat \thetab}_{\rm U}$, ${\hat \thetab}_{\rm C}$ and the tests~$\phi_n$ are such that Assumption~B holds, under a parametric model $\{ {\rm P}_\thetab\n: \thetab \in \Thetab \subset \R^p\}$ that satisfies Assumption~A. Letting~$\lambda (v):=  {\mathbb I}[ v \leq \chi^2_{p-r,1- \alpha}]$, the estimator~${\hat \thetab}_{\rm PTE}$ in \eqref{PTEinvest} rewrites
\begin{equation} 
\label{PTEinvestbis}
{\hat \thetab}_{\rm PTE}:= (1- \lambda(Q\n)) {\hat \thetab}_{\rm U}+  \lambda(Q\n) {\hat \thetab}_{\rm C}.
\end{equation}

%\subsection{Limiting distribution}

\noindent When deriving the asymptotic behavior of~${\hat \thetab}_{\rm PTE}$ under~${\rm P}_\thetab\n$, $\thetab \in \Thetab$, we will discriminate between three cases: (i) $\thetab$ is fixed in the constraint~$\Thetab_0$, (ii) $\thetab=\thetab_n$ belongs to the $\nub_n$-vicinity of the constraint (that is, $\thetab_n= \thetab+ \nub_n \taub_n$, with~$\thetab \in \Thetab_0$ and $(\taub_n)=O(1)$), and (iii) $\thetab$ is fixed outside the constraint~$\Thetab_0$; see Figure~\ref{Asysetup}. 
\vspace{3mm}

\begin{figure}
\begin{center}
\begin{tikzpicture}
 \draw[->] (0,0)--(8,0) node[below right]{$\theta_1$};
  \draw[->] (0,0)--(0,8) node[below left]{$\theta_2$};
    \draw[->] (0,0)--(6,6) node[above right]{$\Thetab_0$};
   \node[circle,fill=black,inner sep=0pt,minimum size=3pt,label=below right:{\scriptsize{(i): $\thetab \in \Thetab_0$}}] (a) at (3,3) {};
   \node[circle,fill=black,inner sep=0pt,minimum size=3pt,label=below right:{\scriptsize{(iii): $\thetab \notin \Thetab_0$}}] (a) at (6.5,2) {};
   \draw[solid] (3,3) circle (1.4);
    \node[circle,fill=black,inner sep=0pt,minimum size=3pt,label=above left:{\scriptsize{(ii): $\thetab+ \nub_n \taub_n$}}] (a) at (2,4) {};
   \draw[-] (3,3)--(2,4); 
\end{tikzpicture}
\end{center}
\caption{Illustration of the various situations where asymptotics are derived on a bivariate parameter $\thetab={\theta_1\choose \theta_2}$, for a constraint of the form~$\Thetab_0={\cal M}(\Umb)$, with $\Umb:={1 \choose 1}$.}
\label{Asysetup}
\end{figure}
\vspace{3mm}

\noindent Our first result shows that, in case~(iii), ${\hat \thetab}_{\rm PTE}$ is asymptotically equivalent to the unconstrained estimator ${\hat \thetab}_{\rm U}$ (see the appendix for a proof).

\newpage
\begin{theorem}
	 \label{Awayresult}
Let Assumptions A and B hold. Fix $\thetab\notin\Thetab_0$ and assume that ${\hat \thetab}_{\rm C}=O_{\rm P}(1)$  under~${\rm P}\n_\thetab$. Then, 
$
\nub_n^{-1}({\hat \thetab}_{\rm PTE}-\thetab)
=
\nub_n^{-1}({\hat \thetab}_{\rm U}-\thetab)+o_{\rm P}(1)
$
 under~${\rm P}\n_\thetab$.
\end{theorem}

We now move to cases (i)--(ii), where we will actually consider parameter sequences of the form $\thetab_n=\thetab+ \nub_n \taub_n \in \Thetab$,  with~$\thetab \in \Thetab_0$ and $(\taub_n)\to\taub$ (note that case~(i) is obtained for~$\taub_n\equiv \mathbf{0}$). We have the following result (see the appendix for a proof).

\begin{theorem}
	\label{Vicinityresult}
Let Assumptions A and B hold and consider sequences of the form $\thetab_n=\thetab+ \nub_n \taub_n \in \Thetab$,  with~$\thetab \in \Thetab_0$ and $(\taub_n)\to\taub$. Conditional on~${\bf D}\n={\bf D}$, $\nub_n^{-1} (\hat\thetab_{\rm PTE}-\thetab_n)$ is, under~${\rm P}\n_{\thetab_n}$, asymptotically normal with mean vector
\begin{eqnarray}
\lefteqn{
\hspace{-9mm}
{\pmb \mu}_{\rm PTE}^{\rm Vic} 
=
(1- \lambda (\| {\bf D} \|^2)) 
\,
\big\{
(\Ab_\thetab \Omegab_\thetab - {\bf I}_p)\taub+ {\bf A}_\thetab \Sigb_\thetab {\bf C}_\thetab\pr ({\bf C}_\thetab \Sigb_\thetab {\bf C}_\thetab\pr)^{-} ({\bf D}- {\bf C}_\thetab \Omegab_\thetab \taub)
\big\} 
}
\nonumber 
\\[1mm] 
& & 
\hspace{-11mm}
+
 \lambda (\| {\bf D} \|^2) 
 \,
\big\{
( \Umb  {\bf B}_\thetab\Omegab_\thetab -{\bf I}_p) \taub+ \Umb{\bf B}_\thetab   \Sigb_\thetab {\bf C}_\thetab \pr ({\bf C}_\thetab \Sigb_\thetab {\bf C}_\thetab\pr)^{-} ({\bf D}- {\bf C}_\thetab \Omegab_\thetab \taub) 
\big\} 
\end{eqnarray}
and covariance matrix
$$
\Gamb_{\rm PTE}^{\rm Vic}
%&=&
%(1- \lambda (\| {\bf D} \|^2))^2 (\Sigb_\thetab^{-1} - \Sigb_\thetab^{-1/2} {\bf P}_{\Umb} \Sigb_\thetab^{-1/2})
%+
% \lambda^2 (\| {\bf D} \|^2) {\bf A}_\Umb  + 2  \lambda (\| {\bf D} \|^2) (1- \lambda (\| {\bf D} \|^2)) {\bf A}_\Umb
%\\[1mm]
=
 (1- \lambda (\| {\bf D} \|^2)) {\bf A}_\thetab(\Sigb_\thetab- {\bf L}_\thetab){\bf A}_\thetab\pr +  \lambda (\| {\bf D} \|^2) \Umb{\bf B}_\thetab (\Sigb_\thetab- {\bf L}_\thetab) {\bf B}_\thetab\pr \Umb\pr,
$$
where we denoted as~${\bf A}^-$ the Moore-Penrose inverse of~${\bf A}$ and where we let~${\bf L}_\thetab:=\Sigb_\thetab {\bf C}_\thetab\pr ( {\bf C}_\thetab \Sigb_\thetab {\bf C}_\thetab\pr)^{-}  {\bf C}_\thetab \Sigb_\thetab$.
\end{theorem}

\noindent Theorem~\ref{Vicinityresult} allows us to obtain an expression for the unconditional asymptotic distribution of $\nub_n^{-1} (\hat\thetab_{\rm PTE}-\thetab_n)$: in the framework of Assumption~(B), the Le Cam third lemma implies that~$\Db_n$ is asymptotically normal with mean vector~$\mub_\Db:={\bf C}_\thetab \Omegab_\thetab \taub$ and covariance matrix~$\Sigb_\Db={\bf C}_\thetab \Sigb_\thetab {\bf C}_\thetab\pr$ under~${\rm P}\n_{\thetab_n}$, so that, under the same sequence of hypotheses, $\nub_n^{-1} (\hat\thetab_{\rm PTE}-\thetab_n)$ converges weakly to a random $p$-vector ${\bf Z}$ with probability density function (pdf)
\begin{equation}\label{limuncond}
{\bf z} \mapsto \int_{\R^p} \phi_{{\pmb \mu}_{\rm PTE}^{\rm Vic}, \Gamb_{\rm PTE}^{\rm Vic}}({ \bf z}) \phi_{\mub_\Db,\Sigb_\Db}(\xb) d\xb,
\end{equation}
where $\phi_{\mub,\Sigb}$ stands for the pdf of the $p$-variate normal distribution with mean vector~$\mub$ and covariance matrix~$\Sigb$. Since the pdf~\eqref{limuncond} does not allow for a simple comparison between~$\hat\thetab_{\rm PTE}$, $\hat\thetab_{\rm U}$ and $\hat\thetab_{\rm C}$, we will base such a comparison on the asymptotic mean square errors (MSEs) of these estimators. 

A general expression for the asymptotic MSEs can be obtained by computing~${\rm E}[{\pmb \mu}_{\rm PTE}^{\rm Vic}]$, ${\rm Var}[{\pmb \mu}_{\rm PTE}^{\rm Vic}]$ and ${\rm E}[\Gamb_{\rm PTE}^{\rm Vic}]$, recalling that, under~${\rm P}\n_{\thetab_n}$, the random vector~$\Db\n$ has asymptotic mean~$\mub_\Db$ and covariance matrix~$\Sigb_\Db$. We now derive these limiting MSEs when PTEs are based on~$Q\n$ in \eqref{optitest} and on asymptotically efficient estimators satisfying~\eqref{optiU}--\eqref{optiC} (limiting MSEs of PTEs based on other estimators can be obtained in the same way). For such estimators and tests, Theorem~\ref{Vicinityresult} yields that, conditional on~${\bf D}\n={\bf D}$,  $\nub_n^{-1} (\hat\thetab_{\rm PTE}-\thetab_n)$ is, under~${\rm P}\n_{\thetab_n}$, asymptotically normal with mean vector
%\begin{eqnarray*}
%{\pmb \mu}_{\rm PTE, eff}^{\rm Vic} 
%&\!\!=\!\!& (1- \lambda (\| {\bf D} \|^2)) (\Gamb_\thetab^{-1/2} {\bf P}^{\perp}_{\Umb}({\bf D}-  \Gamb_\thetab^{1/2} \taub)) 
%\\
%& & 
%\hspace{43mm}  
%+
%\lambda (\| {\bf D} \|^2) (\Gamb_\thetab^{-1/2} {\bf P}_{\Umb}\Gamb_\thetab^{1/2}- {\bf I}_p) \taub 
%\\
%&\!\!=\!\!&  \Gamb_\thetab^{-1/2} {\bf P}^{\perp}_{\Umb} 
%\big\{
%(1- \lambda (\| {\bf D} \|^2)) {\bf D}- \Gamb_\thetab^{1/2} \taub
%\big\}
%,
%\end{eqnarray*}
\begin{equation}
\label{109}	
{\pmb \mu}_{\rm PTE, eff}^{\rm Vic} 
=
  \Gamb_\thetab^{-1/2} {\bf P}^{\perp}_{\Umb,{\rm eff}} 
\big\{
(1- \lambda (\| {\bf D} \|^2)) {\bf D}- \Gamb_\thetab^{1/2} \taub
\big\}
,
\end{equation}
and covariance matrix
%\begin{eqnarray}
%\Gamb_{\rm PTE, eff}^{\rm Vic} 
%&\!\!=\!\!& (1- \lambda (\| {\bf D} \|^2)) (\Gamb_\thetab^{-1}- \Gamb_\thetab^{-1/2} {\bf P}^{\perp}_{\Umb} \Gamb_\thetab^{-1/2})+ \lambda (\| {\bf D} \|^2) \Gamb_\thetab^{-1/2} {\bf P}_{\Umb} \Gamb_\thetab^{-1/2} 
%\nonumber 
%\\
%&\!\!=\!\!& 
%%\Gamb_\thetab^{-1}- \Gamb_\thetab^{-1/2} {\bf P}^{\perp}_{\Umb} \Gamb_\thetab^{-1/2}=
%\Gamb_\thetab^{-1/2} {\bf P}_{\Umb} \Gamb_\thetab^{-1/2}
%\label{110}
%,
%\end{eqnarray}
\begin{equation}
\Gamb_{\rm PTE, eff}^{\rm Vic} 
=
\Gamb_\thetab^{-1/2} {\bf P}_{\Umb,{\rm eff}} \Gamb_\thetab^{-1/2}
\label{110}
,
\end{equation}
with~${\bf P}_{\Umb,{\rm eff}}:=\Gamb_{\thetab}^{1/2}\Umb(\Umb\pr \Gamb_{\thetab} \Umb)^{-1}\Umb\pr \Gamb_{\thetab}^{1/2}$ and~${\bf P}^{\perp}_{\Umb,{\rm eff}}:={\bf I}_p-{\bf P}_{\Umb,{\rm eff}}$. 
%(the last inequality in~(\ref{110}) follows from the identity~${\bf P}^{\perp}_{\Umb}+ {\bf P}_\Umb= {\bf I}_p$). 
We then have the following result (see the appendix for a proof).

\begin{proposition} 
\label{measqprop}
If~${\pmb \mu}_{\rm PTE, eff}^{\rm Vic}$ in~(\ref{109}) is based on a random $p$-vector~$\Db$ that is normal with mean vector~${\bf P}^{\perp}_{\Umb} \Gamb_\thetab^{1/2} \taub$ and covariance matrix~${\bf P}^{\perp}_{\Umb}$, then  
$$
{\rm E}[{\pmb \mu}_{\rm PTE, eff}^{\rm Vic}]= -\gamma_2 \Gamb_\thetab^{-1/2}{\bf P}^{\perp}_{\Umb,{\rm eff}} \Gamb_\thetab^{1/2} \taub
$$ 
and
\begin{eqnarray*}
	\lefteqn{
\hspace{-3mm} 
{\rm Var}[{\pmb \mu}_{\rm PTE, eff}^{\rm Vic}]
= 
(1- \gamma_2)
\Gamb_\thetab^{-1/2} {\bf P}^{\perp}_{\Umb,{\rm eff}}\Gamb_\thetab^{-1/2}
}
\\[2mm]
& & 
\hspace{3mm} 
+ ((1- \gamma_4)-(1- \gamma_2)^2) 
\Gamb_\thetab^{-1/2} {\bf P}^{\perp}_{\Umb,{\rm eff}}\Gamb_\thetab^{1/2} \taub\taub\pr \Gamb_\thetab^{1/2}
{\bf P}^{\perp}_{\Umb,{\rm eff}}\Gamb_\thetab^{-1/2}
,
\end{eqnarray*}
where we let~$\gamma_j:={\rm P}[V_j \leq \chi^2_{p-r,1- \alpha}]$, with~$V_j\sim\chi^2_{p-r+j}(\taub\pr \Gamb_\thetab^{1/2} {\bf P}_{\Umb,{\rm eff}}^\perp \Gamb_\thetab^{1/2} \taub)$ $($throughout, $\chi^2_{\ell}(\eta)$ will stand for the non-central chi-square distribution with~$\ell$ degrees of freedom and with non-centrality parameter~$\eta)$.
\end{proposition}

We define the limiting MSE of~$\hat\thetab_{\rm PTE}$ under~${\rm P}\n_{\thetab_n}$ as
$$
{\rm AMSE}_{\thetab,\taub}(\hat\thetab_{\rm PTE})
:= 
{\rm E}[\Zb \Zb \pr]={\rm Var}[\Zb]+ {\rm E}[\Zb]({\rm E}[\Zb])\pr
,
$$
where~${\bf Z}$ is the weak limit of~$\nub_n^{-1} (\hat\thetab_{\rm PTE}-\thetab_n)$ under~${\rm P}\n_{\thetab_n}$.
Now, since
$
{\rm E}[\Zb]
=
{\rm E}[{\rm E}[\Zb | {\bf D}]]
=
{\rm E}[{\pmb \mu}_{\rm PTE, eff}^{\rm Vic}]
$
and
$
{\rm Var}[\Zb]
=
{\rm E}[{\rm Var}[\Zb | {\bf D}]] + {\rm Var}[{\rm E}[\Zb | {\bf D}]]
=
\Gamb_{\rm PTE, eff}^{\rm Vic} + {\rm Var}[{\pmb \mu}_{\rm PTE, eff}^{\rm Vic}]
$
(note that~${\rm Var}[\Zb | {\bf D}]=\Gamb_{\rm PTE, eff}^{\rm Vic}$ is non-random), Proposition~\ref{measqprop} yields
\begin{eqnarray} 
{\rm AMSE}_{\thetab,\taub}(\hat\thetab_{\rm PTE})
&\!=\!& 
\Gamb_{\rm PTE, eff}^{\rm Vic}+{\rm Var}[{\pmb \mu}_{\rm PTE, eff}^{\rm Vic}]+ ({\rm E}[{\pmb \mu}_{\rm PTE, eff}^{\rm Vic}])({\rm E}[{\pmb \mu}_{\rm PTE, eff}^{\rm Vic}])\pr 
\nonumber 
\\
&\!=\!& 
\Gamb_\thetab^{-1}- \gamma_2 \Gamb_\thetab^{-1/2} {\bf P}^{\perp}_{\Umb,{\rm eff}} \Gamb_\thetab^{-1/2}  
\nonumber  
\\ 
& & 
\hspace{-2mm} 
+ (2 \gamma_2- \gamma_4) \Gamb_\thetab^{-1/2}{\bf P}^{\perp}_{\Umb,{\rm eff}}\Gamb_\thetab^{1/2} \taub \taub\pr \Gamb_\thetab^{1/2}{\bf P}^{\perp}_{\Umb,{\rm eff}} \Gamb_\thetab^{-1/2}
\!
. 
\label{amsePTE}
\end{eqnarray}

To enable proper comparison with the unconstrained and constrained antecedents of~$\hat\thetab_{\rm PTE}$ (namely, the estimators~$\hat\thetab_{\rm U}$ and~$\hat\thetab_{\rm C}$ satisfying~\eqref{optiU} and~\eqref{optiC}, respectively), the following result provides explicit expressions for the limiting MSEs of these estimators (see the appendix for a proof).

\begin{proposition} 
\label{amseUC} 
Let Assumptions A and B hold. Then, under ${\rm P}_{\thetab_n}\n$,
$$
{\rm AMSE}_{\thetab,\taub}(\hat\thetab_{\rm U})=\Gamb_\thetab^{-1}
$$
and 
$$
{\rm AMSE}_{\thetab,\taub}(\hat\thetab_{\rm C})
=
\Gamb_\thetab^{-1/2} {\bf P}_{\Umb,{\rm eff}} \Gamb_\thetab^{-1/2}+ \Gamb_\thetab^{-1/2}{\bf P}_{\Umb,{\rm eff}}^{\perp}\Gamb_\thetab^{1/2} \taub \taub\pr \Gamb_\thetab^{1/2}{\bf P}_{\Umb,{\rm eff}}^{\perp} \Gamb_\thetab^{-1/2}
,
$$
where~$\hat\thetab_{\rm U}$ and $\hat\thetab_{\rm C}$ are estimators satisfying~\eqref{optiU} and~\eqref{optiC}, respectively.
\end{proposition}

\noindent 
It is worthwile to consider some boundary cases. 
%
%(a)
For~$\alpha=1$, we have~$\gamma_2= \gamma_4=0$, so that~${\rm AMSE}_{\thetab,\taub}(\hat\thetab_{\rm PTE})={\rm AMSE}_{\thetab,\taub}(\hat\thetab_{\rm U})$, which is compatible with the fact that~$\hat\thetab_{\rm PTE}=\hat\thetab_{\rm U}$ almost surely when the test~$\phi_n$ is performed at asymptotic level~$\alpha=1$. 
%
%(b)
At the other extreme, for~$\alpha=0$, we rather have~$\gamma_2= \gamma_4=1$, which provides
\begin{eqnarray*}
{\rm AMSE}_{\thetab,\taub}(\hat\thetab_{\rm PTE})
&\!=\!&
\Gamb_\thetab^{-1}-\Gamb_\thetab^{-1/2} {\bf P}^{\perp}_{\Umb,{\rm eff}} \Gamb_\thetab^{-1/2}+\Gamb_\thetab^{-1/2}{\bf P}^{\perp}_{\Umb,{\rm eff}}\Gamb_\thetab^{1/2} \taub \taub\pr \Gamb_\thetab^{1/2}{\bf P}^{\perp}_{\Umb,{\rm eff}} \Gamb_\thetab^{-1/2} 
\\
&\!=\!&
\Gamb_\thetab^{-1/2} {\bf P}_{\Umb,{\rm eff}} \Gamb_\thetab^{-1/2}+\Gamb_\thetab^{-1/2}{\bf P}^{\perp}_{\Umb,{\rm eff}}\Gamb_\thetab^{1/2} \taub \taub\pr \Gamb_\thetab^{1/2}{\bf P}^{\perp}_{\Umb,{\rm eff}} \Gamb_\thetab^{-1/2} 
\\
&\!=\!&
{\rm AMSE}_{\thetab,\taub}(\hat\thetab_{\rm C})
,
\end{eqnarray*} 
in agreement with the fact that~$\hat\thetab_{\rm PTE}=\hat\thetab_{\rm C}$ almost surely when the test~$\phi_n$ is performed at asymptotic level~$\alpha=0$.  

\begin{figure}[!h]
%\vspace*{-0.5cm}
\centering
  \includegraphics[width=0.9\textwidth]{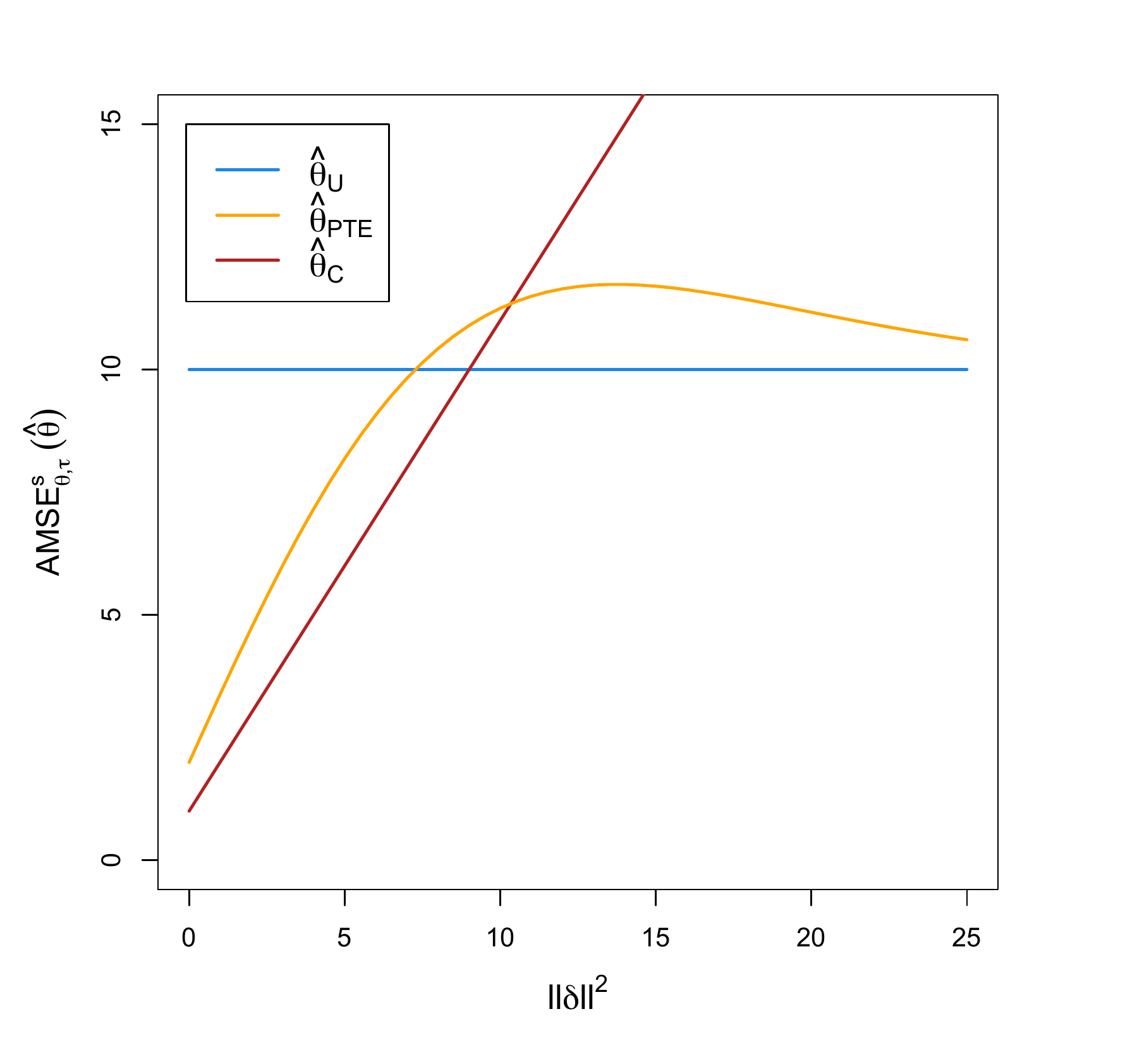}
\vspace*{-5mm}
  \caption{Plots of~${\rm AMSE}^{\rm s}_{\thetab,\taub}(\hat\thetab_{\rm U})$, ${\rm AMSE}^{\rm s}_{\thetab,\taub}(\hat\thetab_{\rm C})$ and~${\rm AMSE}^{\rm s}_{\thetab,\taub}(\hat\thetab_{\rm PTE})$ as functions of~$\| \deltab \|^2$, for $p=10$, $r=1$ and $\alpha=.05$.}
\label{Fig2}
\end{figure}

To conclude this section, we offer a comparison between~${\rm AMSE}_{\thetab,\taub}(\hat\thetab_{\rm PTE})$, ${\rm AMSE}_{\thetab,\taub}(\hat\thetab_{\rm U})$, and~${\rm AMSE}_{\thetab,\taub}(\hat\thetab_{\rm C})$. These limiting MSEs being matrix-valued, it is needed to base this comparison on a scalar summary, such as, e.g., their trace. In the present case, where the unconstrained estimator satisfies~${\rm AMSE}_{\thetab,\taub}(\hat\thetab_{\rm U})=\Gamb_\thetab^{-1}$, it is natural to measure the asymptotic performance of an estimator~$\hat{\thetab}$ through the equivalent scalar quantity
$$
{\rm AMSE}^{\rm s}_{\thetab,\taub}(\hat\thetab)
:=
{\rm tr}[\Gamb_\thetab^{1/2}({\rm AMSE}_{\thetab,\taub}(\hat\thetab)) \Gamb_\thetab^{1/2}]
,
$$
which, for~$\hat\thetab_{\rm U}$, will provide the ``normalized" perfomance~${\rm AMSE}^{\rm s}_{\thetab,\taub}(\hat\thetab_{\rm U})=p$ (see Proposition~\ref{amseUC}), that does not depend on the value of~$\thetab$ at which the contiguous alternatives~$\thetab_n=\thetab+ \nub_n \taub_n$ are localized.    
%\begin{equation} \label{Frofro1}
%<\Gamb_\thetab, {\rm AMSE}_{\thetab,\taub}(\hat\thetab_{\rm U})>= {\rm tr}(\Gamb_\thetab \Gamb_\thetab^{-1})=p.
%\end{equation}
%Since ${\rm AMSE}_{\thetab,\taub}(\hat\thetab_{\rm U})$ does not depend on $\taub$, \eqref{Frofro1} holds for any bounded sequence $\taub_n$ and therefore in particular when $\taub_n \equiv {\bf 0}$. 
Proposition~\ref{amseUC} also entails that 
$$
{\rm AMSE}^{\rm s}_{\thetab,\taub}(\hat\thetab_{\rm C}) 
=
{\rm tr}[{\bf P}_{\Umb,{\rm eff}}]+{\rm tr}[{\bf P}_{\Umb,{\rm eff}}^{\perp}\Gamb_\thetab^{1/2} \taub \taub\pr \Gamb_\thetab^{1/2}{\bf P}_{\Umb,{\rm eff}}^{\perp}] 
=
r+ \| \deltab\|^2
, 
$$
with~$\deltab:={\bf P}_{\Umb,{\rm eff}}^{\perp}\Gamb_\thetab^{1/2} \taub$. Note that, at~$\taub=0$, this shows that~$
{\rm AMSE}^{\rm s}_{\thetab,\taub}(\hat\thetab_{\rm C})=r<p={\rm AMSE}^{\rm s}_{\thetab,\taub}(\hat\thetab_{\rm U})$, which confirms the intuition that~$\hat{\thetab}_{\rm C}$ dominates~$\hat{\thetab}_{\rm U}$ when the true parameter value belongs to~$\Thetab_0$. Now, it easily follows from~(\ref{amsePTE}) that
$$
{\rm AMSE}^{\rm s}_{\thetab,\taub}(\hat\thetab_{\rm PTE}) 
=
p-\gamma_2 (p-r)+(2 \gamma_2- \gamma_4) \| \deltab \|^2
,
$$
where $\gamma_j={\rm P}[V_j \leq \chi^2_{p-r,1- \alpha}]$, with $V_j \sim \chi^2_{p-r+j}(\| \deltab \|^2)$. Figure~\ref{Fig2}  plots, for $p=10$, $r=1$ and $\alpha=.05$, the quantities~${\rm AMSE}^{\rm s}_{\thetab,\taub}(\hat\thetab_{\rm U})$, ${\rm AMSE}^{\rm s}_{\thetab,\taub}(\hat\thetab_{\rm C})$ and~${\rm AMSE}^{\rm s}_{\thetab,\taub}(\hat\thetab_{\rm PTE})$ as functions of~$\| \deltab \|^2$. The figure reveals that, under~${\rm P}_\thetab\n$ with $\thetab \in \Thetab_0$ (which corresponds to~$\deltab=0$), the constrained estimator~$\hat{\thetab}_{\rm C}$ has the best performance, as expected. The PTE performs better than~$\hat{\thetab}_{\rm U}$ in the vicinity of the constraint ($\|\deltab\|$ small to moderate) and it is asymptotically equivalent to~$\hat{\thetab}_{\rm U}$ far from the constraint ($\|\deltab\|$ large).

%When $\taub={\bf 0}$, we have that ${\rm tr}(\Gamb_\thetab {\rm AMSE}_\thetab(\hat\thetab_{\rm U}))=p$, ${\rm tr}(\Gamb_\thetab {\rm AMSE}_\thetab(\hat\thetab_{\rm C}))=r$ and  ${\rm tr}(\Gamb_\thetab {\rm AMSE}_\thetab(\hat\thetab_{\rm PTE}))=p-\gamma_2 (p-r)$, where $\gamma_2={\rm P}[V \leq \chi^2_{p-r,1- \alpha}]$ with $V \sim \chi^2_{p-r+2}$. 
%When $\taub \neq {\bf 0}$, putting $\deltab:={\bf P}^{\perp}_{\Umb}\Gamb_\thetab^{1/2} \taub$, we have that ${\rm tr}(\Gamb_\thetab {\rm AMSE}_\thetab(\hat\thetab_{\rm U}))=p$, ${\rm tr}(\Gamb_\thetab {\rm AMSE}_\thetab(\hat\thetab_{\rm C}))=r+ \| \deltab \|^2$ and ${\rm tr}(\Gamb_\thetab{\rm AMSE}_\thetab(\hat\thetab_{\rm PTE}))=p-\gamma_2 (p-r)+(2 \gamma_2- \gamma_4) \| \deltab \|^2$,
%where $\gamma_j={\rm P}[V_j \leq \chi^2_{p-r,1- \alpha}]$ with $V_j \sim \chi^2_{p-r+j}(\| \deltab \|^2)$. In Figure \ref{Fig2} we provide a plot of ${\rm tr}(\Gamb_\thetab{\rm AMSE}_\thetab(\hat\thetab))$ for $\hat\thetab=\hat\thetab_{\rm U}, \hat\thetab=\hat\thetab_{\rm C}$ and $\hat\thetab=\hat\thetab_{\rm PTE}$ as a function of $\| \deltab \|^2$ for $p=3$, $r=1$ and $\alpha=.05$. Inspection of Figure 2\ref{Fig2} reveals that under ${\rm P}_\thetab$ with $\thetab \in {\cal M}(\Umb)$, $\hat{\thetab}_{\rm C}$ has the better performance as expected. The preliminary test estimator performs better that $\hat{\thetab}_{\rm U}$ in the vicinity of the constraint. It is asymptotically equivalent to $\hat{\thetab}_{\rm U}$ away from the constraint.
\vspace{4mm}

\vspace{4mm}

\newpage

\setcounter{section}{4}
\noindent {\bf 4. Two specific applications}
\vspace{4mm}

\noindent In this section, we illustrate the general results obtained above on two particular cases. First, we consider preliminary test estimation in the simple linear regression model and show that we recover for this model and for the considered estimation problem the classical results of \cite{SA06} (Section~4.1). Then, we consider the joint estimation of~$m$ covariance matrices $\Sigb_1, \ldots, \Sigb_m$ in a context where it is suspected that these covariance matrices might be equal, might be proportional, or might share a common ``scale" (Section~4.2).

\subsection{Simple linear regression}

Consider the simple linear regression model
\begin{equation} 
\label{linearmodel}
\Yb= \rho {\bf 1}_n+ \beta \xb+ {\pmb \epsilon}
,
\end{equation}
where $\Yb=(Y_1, \ldots, Y_n)\pr$ is a response vector, $\xb=(x_1, \ldots, x_n)\pr$ is a vector of non-random covariates, and where the error vector~${\pmb \epsilon}=(\epsilon_1, \ldots, \epsilon_n)\pr$ is multinormal with mean zero and covariance matrix~$\sigma^2 {\bf I}_n$, for some $\sigma^2>0$.  This is the classical simple linear model with intercept~$\rho$, slope~$\beta$, and Gaussian homoscedastic errors with variance~$\sigma^2$. Throughout, we consider the parameter~{$\thetab:=(\rho, \beta)'$}, as we will assume that~$\sigma^2$ is known (this is actually no restriction, since the block-diagonality of the Fisher information matrix in this model entails that replacing~$\sigma^2$ with a root-$n$ consistent estimator will have no asymptotic cost, so that all results we obtain below extend to the case where~$\sigma^2$ would remain an unspecified nuisance). One can easily show that this model is ULAN, with a central sequence~$\Deltab_\thetab\n$ that, under ${\rm P}_\thetab\n$, is asymptotically normal with mean zero and covariance matrix
$$
\Gamb_\thetab  
=
 \frac{1}{\sigma^2}  
 \Bigg(
\begin{array}{cc}
1  & \bar{x}_0 \\[-2mm]
\bar{x}_0 & s_0 + \bar{x}_0^2
\end{array}
\Bigg)
,
$$
where $\bar{x}_0:= \lim_{\ny} n^{-1} \sum_{i=1}^n x_i$ and $s_0:= \lim_{\ny} s_x\n$ with $s_x\n:=n^{-1} \xb\pr \xb - n^{-2}({\bf 1}_n\pr \xb)^2$; of course, we tacitly assume that these limits exist and are finite. We consider here preliminary test estimation of~$\thetab$ when it is suspected that~$\beta=\beta_0$ for some given~$\beta_0$. In the context, the classical, unconstrained, estimator of~$\thetab$ is the maximum likelihood estimator 
$$
\hat{\thetab}_{\rm U}
:=
{\hat{\rho} \choose \hat{\beta}}
:=
{n^{-1}({\bf 1}_n\pr \Yb- \hat{\beta} {\bf 1}_n\pr \xb)
 \choose 
(\xb\pr \Yb- n^{-1} \xb\pr {\bf 1}_n {\bf 1}_n\pr \Yb)/ ns_x\n}
,
$$
whereas the natural constrained estimator would be~$\hat{\thetab}_{\rm C}:={\tilde{\rho} \choose \beta_0}$, with~$\tilde{\rho}:= n^{-1}({\bf 1}_n\pr \Yb- {\beta}_0 {\bf 1}_n\pr \xb)$. Since the locally asymptotically optimal test for~$\mathcal{H}_0:\beta=\beta_0$ against~$\mathcal{H}_1:\beta\neq \beta_0$ rejects the null hypothesis at asymptotic level~$\alpha$ when
$$
Q\n
:=
\frac{n(\hat{\beta}- \beta_0)^2 s_x\n}{\sigma^2}
>
 \chi^2_{1, 1- \alpha}
,
$$
the resulting PTE is given by 
$$
\hat{\thetab}_{\rm PTE}
=
{ \hat{\rho}_{\rm PTE} \choose \hat{\beta}_{\rm PTE}}
:=
{\mathbb I}[Q\n >\chi^2_{1, 1- \alpha}] 
\hat{\thetab}_{\rm U}
+ {\mathbb I}[Q\n \leq \chi^2_{1, 1- \alpha}]
\hat{\thetab}_{\rm C}
.
$$
Letting~$\thetab_0={\rho \choose \beta_0}$ be an arbitrary value of the parameter of interest corresponding to the constraint, the null hypothesis can be written as ${\cal H}_0: \thetab \in {\thetab_0} + {\cal M}(\Umb)$, with~$\Umb:={1 \choose 0}$. Since 
$$
\Gamb_{\thetab}^{-1} = \sigma^2
\Bigg(
\begin{array}{cc}
1 + \frac{\bar{x}_0^2}{s_0} & -\frac{\bar{x}_0}{s_0} \\[-3mm]
-\frac{\bar{x}_0}{s_0} & \frac{1}{s_0}
\end{array}
\Bigg) 
\quad 
{\rm and} 
\quad
\Gamb_\thetab^{-1/2} {\bf P}^{\perp}_{\Umb} \Gamb_\thetab^{1/2}
=
\Bigg(
\begin{array}{cc}
 0 & -\bar{x}_0 \\[-2mm]
  0 & 1 
 \end{array}
 \!
\Bigg) 
,
$$
it follows from \eqref{amsePTE} that, under~${\rm P}\n_{\thetab_0+ n^{-1/2}\taub}$, with~$\taub={0\choose \delta}$,  the MSE quantity~${\rm AMSE}_{\thetab,\taub}(\hat{\thetab}_{\rm PTE})$ is here given by
$$
  \Bigg(
\begin{array}{cc} 
\sigma^2(1+\frac{\bar{x}^2_0}{s_0}- \frac{\gamma_2 \bar{x}_0^2}{s_0})+ (2 \gamma_2- \gamma_4)  \bar{x}_0^2 \delta^2  
&
  \frac{\sigma^2(\gamma_2-1)\bar{x}_0}{s_0}- (2 \gamma_2- \gamma_4)  \bar{x}_0 \delta^2 
  \\[-2mm]  
  \frac{\sigma^2(\gamma_2-1)\bar{x}_0}{s_0}- (2 \gamma_2- \gamma_4)  \bar{x}_0 \delta^2 
  & 
  \frac{ \sigma^2(1-\gamma_2)}{s_0}+ (2 \gamma_2- \gamma_4)  \delta^2 
 \end{array}
 \!
\Bigg) 
,
$$
 where the $\gamma_j$'s are computed with $p=2$ and $r=1$. This is in perfect agreement with the result 
%\begin{eqnarray}
%{\rm AMSE}((\hat{\theta}_{\rm PTE}, \hat{\beta}_{\rm PTE})\pr) &=& \sigma^2  (1- {\bf p}_3)  \begin{pmatrix} 1+\frac{\bar{x}^2_0}{Q^2} & - \frac{\bar{x}_0}{Q^0} \\ - \frac{\bar{x}_0}{Q^0} & \frac{1}{Q^2} \end{pmatrix} 
% + \sigma^2 {\bf p}_3 \begin{pmatrix} 1 & 0 \\ 0 & 0 \end{pmatrix} 
%  \nonumber \\[1mm]
% &+&
%  (2  {\bf p}_3 -  {\bf p}_5)  \begin{pmatrix}  \bar{x}^2_0 \delta^2 & - \bar{x}_0 \delta^2 \\ - \bar{x}_0 \delta^2 & \delta^2 \end{pmatrix},
%  \end{eqnarray}
in Theorem~4, \mbox{p.p.}~94--96 in \cite{SA06}.

\subsection{Multisample estimation of covariance matrices}
\vspace{4mm}

\noindent Consider $m(\geq 2)$ mutually independent samples of random $k$-vectors \linebreak $\Xb_{i1}, \ldots, \Xb_{in_i}$, $i=1, \ldots,m$, with respective sample sizes $n_1, \ldots, n_m$, such that, for any~$i$, the~$\Xb_{ij}$'s form a random sample from the   multinormal distribution with mean vector~${\bf 0}$ and (invertible) covariance matrix~$\Sigb_i$ (all results below extend to the case where observations in the~$i$th sample would have a common, unspecified, mean~${\pmb\mu}_{i}$, $i=1,\ldots,n$, due to the block-diagonality of the Fisher information matrix for location and scatter in elliptical models; see, e.g., \citealp{HP06}). In the sequel, we decompose the covariance matrices into~$\Sigb_i=\sigma^2_i \Vb_i$, where~$\sigma_i^2:= ({\rm det}\,\Sigb_i)^{1/k}$ is their ``scale"  and ${\Vb}_i:= \Sigb_i/ ({\rm det}\,\Sigb_i)^{1/k}$ is their ``shape". Under the only assumption that~$\lambda_i:=\lambda_i\n:=n_i/n:=n_i/(\sum_{\ell=1}^m n_\ell) \to \lambda_i \in (0,1)$ (to make the notation lighter, we will not stress the dependence in~$n$ in many quantities below), it follows from \cite{HP09} that the sequence of Gaussian models indexed by 
$$
\thetab
:=
\left(\sigma_1^{2},\ldots,\sigma_m^{2},({\vech }{\bf V }_1)\pr,\ldots,({\vech}
{\bf V }_m)\pr\right)\pr,
$$
where~$\vech \Vb(\in\R^{d_k}$, with~$d_k:=k(k+1)/2-1)$ stands for the vector obtained by depriving~${\rm vech} \Vb$ of its first entry~$\Vb_{11}$, is ULAN in the sense of Assumption~A. To describe the corresponding central sequence and Fisher information matrix, we need the following notation. 

Denoting as ${\bf e}_r$ the $r$th vector of the canonical basis of~$\R^k$,  we let ${\bf K}_{k}:= \sum_{r,s=1}^{k}({\bf e}_r{\bf e}_s\pr)\otimes ({\bf e}_s{\bf e}_r\pr)$ be the $k^2\times k^2$ {\em commutation
matrix}, put ${\bf J}_k: = ({\rm vec}\,{\bf I}_k) ({\rm
vec}\,{\bf I}_k)\pr$, and define~${\bf M}_k({\bf V})$ \label{definmk} as the $(d_k \times k^2)$ matrix  such that $({\bf M}_k({\bf V}))\pr (\vech {\bf v}) = {\rm vec}\,{\bf v}$ for any symmetric~$k\times k$ matrix~${\bf v}$ such that~${\rm tr}[{\bf V}^{-1}{\bf v}] = 0$. We further put
\vspace{-2mm}
$$
{\bf H}_k({\bf V}):=\frac{1}{4} \,{\bf M}_k({\bf V})  \left({\bf I}_{k^2} + {\bf K}_k  \right) \left( {\bf V}^{\otimes 2}\right)^{-1}  ({\bf M}_k({\bf V}))\pr.
$$
 Then, letting $\Sb_i:=n_i^{-1} \sum_{j=1}^{n_i} \Xb_{ij} \Xb_{ij}\pr$ be the empirical covariance matrix in sample $i$, the central sequence is~$\Deltab_\thetab=(\Delta_\thetab^{\I,1}, \ldots, \Delta_\thetab^{\I,m}, \Deltab_\thetab^{\II,1}, \ldots, \Deltab_\thetab^{\II,m})$, where, for~$i=1,\ldots , m$, we wrote
$$
\Delta^{\I,i} _{{\thetab}} 
:= 
\frac{ \sqrt{n_i}  }{2\sigma_i^{2}} 
\,
{\rm tr}\big[\sigma_i^{-2} \Vb_i^{-1} (\Sb_i- \sigma_i^2 \Vb_i)\big]
, 
\vspace{-1mm}
\ \
\Deltab^{\II,i} _{{\pmb\theta}} 
:= 
\frac{ \sqrt{n_i}  }{2 \sigma_i^{2}}  
\,
{\bf M}_k({\bf V}_i)
\!\left( {\bf V}_i^{\otimes 2}
\right)^{-1}
\!
({\rm vec}\,{\Sb_i})
, 
$$
whereas the (full-rank and block-diagonal) information matrix takes the form
$\Gamb  _{{\pmb\theta}} 
:=
{\rm diag}(\Gamb^{\I}_{{\pmb\theta}} , \Gamb^{\II}_{{\pmb\theta}})
:=
{\rm diag}(\frac{k}{2}\sigbu ^{-4} , {\bf H}_k(\underline{\bf V}))
$,
with~$\sigbu :={\rm diag}(\sigma_1,\ldots,\sigma_m)$ and ${\bf H}_k(\underline{\bf V}) :={\rm diag}({\bf H}_k({\bf
V}_1),\ldots, {\bf H}_k({\bf V}_m))$. 
The corresponding contiguity rate $\nub_n$ in Assumption~A is given by 
$\nub_n= n^{-1/2} {\bf r}_n$, where 
$$
{\bf r}_n
:=
{\rm diag}\big(
\lambda_1^{-1/2}, \ldots, \lambda_m^{-1/2}
,
\lambda_1^{-1/2} {\bf I}_{d_k}, \ldots, \lambda_m^{-1/2} {\bf I}_{d_k}
 \big)
.
$$

We consider here estimation of~$\Sigb_1, \ldots, \Sigb_m$ or, equivalently, estimation of~$\thetab$. An advantage of the $\thetab$-parametrization is that it allows the construction of various PTEs: one may suspect, e.g., \emph{scale homogeneity} ${\cal H}_0^{\rm scale}:\sigma^2_1= \ldots= \sigma^2_m$, \emph{shape homogeneity} ${\cal H}_0^{\rm shape}: \Vb_1= \ldots= \Vb_m$, or full \emph{covariance homogeneity} ${\cal H}_0^{\rm cov}:\sigma^2_1 \Vb_1= \ldots= \sigma_m^2 \Vb_m$, that is, ${\cal H}_0^{\rm cov}:\Sigb_1= \ldots= \Sigb_m$. An asymptotically efficient unconstrained estimator in this Gaussian model is given by
\begin{equation}
\label{unconstrained}
{\hat \thetab}_{\rm U}
:=
\left(
 ({\rm det}\,\Sb_1)^{1/k},\ldots, ({\rm det}\,\Sb_m)^{1/k}, \frac{({\vech }{\Sb }_1)\pr}{({\rm det}\,\Sb_1)^{1/k}},\ldots,\frac{({\vech} {\Sb }_m)\pr}{({\rm det}\,\Sb_m)^{1/k}} 
\right)\pr
,
\end{equation}
whereas, writing~$\Sb:= n^{-1} \sum_{i=1}^m \sum_{j=1}^{n_i} \Xb_{ij} \Xb_{ij}\pr$ for the pooled covariance matrix estimator,  asymptotically efficient constrained estimators, for the three constraints~${\cal H}_0^{\rm scale}$, ${\cal H}_0^{\rm shape}$ and ${\cal H}_0^{\rm cov}$ above, are given by
\begin{equation}
\label{constrained1}
{\hat \thetab}_{\rm C}^{\rm scale}
:=
\left( ({\rm det}\,\Sb)^{1/k} {\bf 1}_m\pr ,\frac{({\vech }{\Sb }_1)\pr}{({\rm det}\,\Sb_1)^{1/k}},\ldots,\frac{({\vech}
{\Sb }_m)\pr}{({\rm det}\,\Sb_m)^{1/k}} \right)\pr,
\end{equation}
\begin{equation}
\label{constrained2}
{\hat \thetab}_{\rm C}^{\rm shape}
:=
\left( ({\rm det}\,\Sb_1)^{1/k},\ldots, ({\rm det}\,\Sb_m)^{1/k}, {\bf 1}_m\pr \otimes \frac{({\vech }{\Sb })\pr}{({\rm det}\,\Sb)^{1/k}} \right)\pr
\end{equation}
and
\begin{equation}
\label{constrained3}
{\hat \thetab}_{\rm C}^{\rm cov}
:=
\left( ({\rm det}\,\Sb)^{1/k} {\bf 1}_m\pr ,  {\bf 1}_m\pr \otimes \frac{({\vech }{\Sb })\pr}{({\rm det}\,\Sb)^{1/k}} \right)\pr
,
\end{equation}
 respectively. The three hypotheses~${\cal H}_0^{\rm scale}$, ${\cal H}_0^{\rm shape}$ and ${\cal H}_0^{\rm cov}$ impose linear restrictions on~$\thetab$, hence can be written as 
 $$
 {\cal H}_0^{\rm scale}: \thetab \in {\cal M}(\Umb_{\rm scale}),
 \quad
 {\cal H}_0^{\rm shape}: \thetab \in {\cal M}(\Umb_{\rm shape}),
 \quad
 {\cal H}_0^{\rm cov}: \thetab \in {\cal M}(\Umb_{\rm cov}),
 \quad
 $$
(more specifically,~$\Umb_{\rm scale}:={\rm diag}({\bf 1}_m, {\bf I}_{md_k})$, $\Umb_{\rm shape}:={\rm diag}({\bf I}_m, {\bf 1}_m \otimes {\bf I}_{d_k})$ and $\Umb_{\rm cov}:={\rm diag}({\bf 1}_m, {\bf 1}_m \otimes {\bf I}_{d_k})$). Now, if the~$p\times r$ matrix~$\Umb$ stands for either of~$\Umb_{\rm scale}$, $\Umb_{\rm shape}$ or $\Umb_{\rm cov}$ (of course, each constraint matrix has its own~$r$), the locally asymptotically most stringent test~$\phi \n_{\Umb}$ for ${\cal H}_0: \thetab \in {\cal M}(\Umb)$ rejects the null hypothesis at asymptotic level~$\alpha$ when
\begin{eqnarray} 
Q\n_{\thetab,\Umb}
& \! := \! &
\Deltab_{\thetab}\pr
\,
\Big[
\Gamb ^{-1}_{\thetab} 
- 
({\bf r}\n)^{-1}{\pmb\Upsilon}
  (\Umb\pr ({\bf r}\n)^{-1} \Gamb_{\thetab} ({\bf r}\n)^{-1} \Umb)^{-1}
   \Umb\pr  ({\bf r}\n)^{-1}  
  \Big]
  \Deltab_{\thetab}
  \nonumber
\\[2mm]
& \! > \! & 
\chi^2_{m(d_k+1)-r,1-\alpha}
.
\label{testgen}
\end{eqnarray}
This allows us  to consider the PTEs
$$
\hat{\thetab}_{\rm PTE}^{\rm scale}:={\mathbb I}[\phi\n_{\Umb_{\rm scale}}=1]{\hat \thetab}_{\rm U} + {\mathbb I}[\phi\n_{\Umb_{\rm scale}}=0]{\hat \thetab}_{\rm C}^{\rm scale}
,
$$ 
$$
\hat{\thetab}_{\rm PTE}^{\rm shape}:={\mathbb I}[\phi\n_{\Umb_{\rm shape}}=1]{\hat \thetab}_{\rm U} + {\mathbb I}[\phi\n_{\Umb_{\rm shape}}=0]{\hat \thetab}_{\rm C}^{\rm shape}
$$ 
and 
$$
\hat{\thetab}_{\rm PTE}^{\rm cov}:={\mathbb I}[\phi\n_{\Umb_{\rm cov}}=1]{\hat \thetab}_{\rm U} + {\mathbb I}[\phi\n_{\Umb_{\rm cov}}=0]{\hat \thetab}_{\rm C}^{\rm cov}
.
$$ 

\begin{figure}[!h]
%\vspace*{-0.5cm}
\centering
  \includegraphics[width=\textwidth,height=60mm]{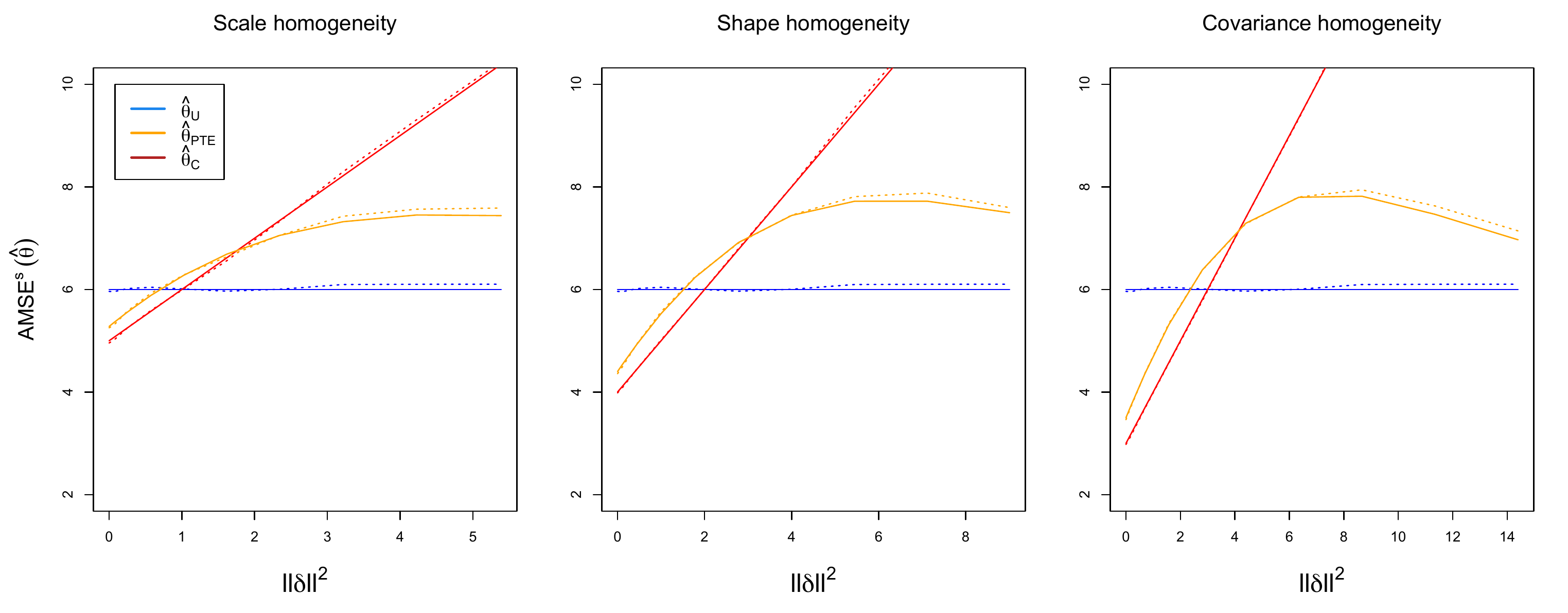}
\vspace*{-9mm}
  \caption{Plots of the empirical performance measures ${\rm tr}[\Gamb_\thetab \hat{\rm AMSE}_\thetab(\hat\thetab)]$ (dotted lines) and of their asymptotic counterparts~${\rm tr}[\Gamb_\thetab{\rm AMSE}_\thetab(\hat\thetab)]$ (solid lines), as functions of~$\| \deltab \|^2$ (which measures distance to the constraint), of the constrained estimators~$\hat{\thetab}_{\rm C}$, unconstrained estimators~$\hat{\thetab}_{\rm U}$ and PTE estimators~$\hat{\thetab}_{\rm PTE}$ associated with the constraints of \emph{scale homogeneity} (left), \emph{shape homogeneity} (middle), and \emph{covariance homogeneity} (right). For the PTEs, all tests are performed at asymptotic level~$\alpha=.05$; see Section~4 for details.}
\label{Fig3}
\end{figure}

To compare these PTEs with their unconstrained and constrained antecedents, we performed the following Monte Carlo exercise, that focuses on the case~$m=2$,~$k=2$ and~$n_1=n_2=20,\!000$. We generated independently~$M=10,\!000$ samples of mutually independent observations $(\Xb_1,\ldots,\Xb_{n1},{\Yb}_1(\ell),\ldots,{\Yb}_{n_2}(\ell))$, $\ell=0, \ldots, 9$, where the~$\Xb_i$'s are ${\cal N}({\bf 0}, {\bf I}_k)$ and the~$\Yb_{i, \ell}$'s are ${\cal N}({\bf 0}, \Sigb_\ell)$, with
$$
\Sigb_\ell
:=
e^{\ell/400}\, 
\frac{{\bf I}_2 +  \ell n^{-1/2} ({\bf e}_2 {\bf e}_1\pr+{\bf e}_1 {\bf e}_2\pr)}{{{\rm det}({\bf I}_2 +  \ell n^{-1/2} ({\bf e}_2 {\bf e}_1\pr+{\bf e}_1 {\bf e}_2\pr))}}
\cdot
$$
For~$\ell=0$, the samples $\Xb_1, \ldots, \Xb_{n_1}$ and $\Yb_1(\ell),\ldots,\Yb_n(\ell)$ share the same underlying covariance matrix~${\bf I}_p$, hence also the same scales and shapes, whereas~$\ell=1, \ldots, 9$ provide  increasingly distinct scales and shapes. In other words, the constraints above are met for~$\ell=0$ and are increasingly violated for~$\ell=1,\ldots,9$. For every considered estimator~$\hat{\thetab}$ of the resulting true parameter value~$\thetab$, we measure the performance of~$\hat{\thetab}$ through 
$$
\hat{\rm AMSE}_\thetab(\hat{\thetab})
:= 
\frac{1}{M}
\sum_{m=1}^M 
(\hat{\rm AMSE}_\thetab(\hat{\thetab}))_{m}
:=
\frac{1}{M}
\sum_{m=1}^M 
n (\hat{\thetab}^{(m)}- \thetab)(\hat{\thetab}^{(m)}- \thetab)\pr
,
$$
where $\hat{\thetab}^{(m)}$ is an estimator computed in the $m$th replication ($m=1, \ldots,M$), or rather, parallel to what we did in Section~3, through the scalar quantity~$\hat{\rm AMSE}^{\rm s}_\thetab(\hat{\thetab}):={\rm tr}[\Gamb_\thetab \hat{\rm AMSE}_\thetab(\hat{\thetab})]$. Figure~\ref{Fig3} then plots~$\hat{\rm AMSE}^{\rm s}(\hat{\thetab})$ for the PTEs~$\hat{\thetab}_{\rm PTE}^{\rm scale}$, $\hat{\thetab}_{\rm PTE}^{\rm shape}$ and~$\hat{\thetab}_{\rm PTE}^{\rm cov}$ (the corresponding tests are all performed at asymptotic level~$\alpha=.05$), as well as their constrained and unconstrained antecedents~${\hat \thetab}_{\rm C}^{\rm scale}$, ${\hat \thetab}_{\rm C}^{\rm shape}$, ${\hat \thetab}_{\rm C}^{\rm cov}$ and~${\hat \thetab}_{\rm U}$. To match what was done in Figure~\ref{Fig2}, these quantities are not plotted as functions of~$\ell$, but rather as functions of the induced~$\| \deltab \|^2$. The figure also provides the corresponding asymptotic performance measures~${\rm AMSE}^{\rm s}_{\thetab,\taub}(\hat{\thetab})$ resulting from the general expression obtained in Section~3. Clearly, the results show that that these empirical and theoretical performance measures are in a perfect agreement.  

%Note that it directly follows from the discussion below Proposition \ref{amseUC} that 
%${\rm tr}(\Gamb_\thetab {\rm AMSE}_\thetab(\hat\thetab_{\rm U}))=m(d_k+1)$, ${\rm tr}(\Gamb_\thetab {\rm AMSE}_\thetab(\hat\thetab_{\rm C (scale)}))=md_k+1$, ${\rm tr}(\Gamb_\thetab {\rm AMSE}_\thetab(\hat\thetab_{\rm PTE (scale)}))=...$, ${\rm tr}(\Gamb_\thetab {\rm AMSE}_\thetab(\hat\thetab_{\rm C (shape)}))=$, ${\rm tr}(\Gamb_\thetab {\rm AMSE}_\thetab(\hat\thetab_{\rm PTE (shape)}))=$, ${\rm tr}(\Gamb_\thetab {\rm AMSE}_\thetab(\hat\thetab_{\rm C (cov)}))=$ and ${\rm tr}(\Gamb_\thetab {\rm AMSE}_\thetab(\hat\thetab_{\rm PTE (cov)}))=$
%with an interesting preliminary test estimator defined as
% \begin{eqnarray}
% {\hat \thetab}_{\rm PTE} &:= &{\mathbb I}[\phi\n_{\Umb_{\rm shape}}=1, \phi\n_{\Umb_{\rm scale}}=1]{\hat \thetab}_{\rm U} + {\mathbb I}[\phi\n_{\Umb_{\rm shape}}=1, \phi\n_{\Umb_{\rm scale}}=0]{\hat \thetab}_{\rm C(scale)} \nonumber \\
% & & + {\mathbb I}[\phi\n_{\Umb_{\rm shape}}=0, \phi\n_{\Umb_{\rm scale}}=1]{\hat \thetab}_{\rm C (shape)}+ {\mathbb I}[\phi\n_{\Umb_{\rm shape}}=0, \phi\n_{\Umb_{\rm scale}}=0]{\hat \thetab}_{\rm C(cov)} \nonumber
 %\end{eqnarray}

\vskip 14pt
\noindent {\large\bf Acknowledgements}
\vskip 14pt

\noindent Davy Paindaveine's research is supported by a research fellowship from the Francqui Foundation and by the Program of Concerted Research Actions (ARC) of the Universit\'{e} libre de Bruxelles. Thomas Verdebout's research is supported by the ARC Program of the Universit\'{e} libre de Bruxelles and by the Cr\'{e}dit de Recherche  J.0134.18 of the FNRS (Fonds National pour la Recherche Scientifique), Communaut\'{e} Fran\c{c}aise de Belgique.
\par
%%%%%%%%%%%%%%%%%%%%%%%%%%%%%%%%%%%%%%%%%%%%%%%%%%%%%%%%%%%%%%%%%%%%%%%%%%%%%%%%%%%%%%%%%%%%%%%%%%%%%%%%%%%%%%%%%%%%%%%%%%%%
\vskip 14pt
\noindent {\large\bf Appendix: Proofs}

\appendix
\setcounter{section}{5}
%\setcounter{equation}{0} %-1
%\noindent {\bf 5. Proofs}
\vspace{4mm}

\noindent In this appendix, we collect the proofs of the various results. \vspace{2mm}

\noindent {\bf Proof of Theorem~\ref{Awayresult}}. First note that
\begin{eqnarray}
\nub_n^{-1} ({\hat \thetab}_{\rm PTE}-\thetab)
&\!\!=\!\!&
 \lambda(Q\n) \nub_n^{-1} ({\hat \thetab}_{\rm C}-\thetab) 
+
(1-\lambda(Q\n)) \nub_n^{-1} ({\hat \thetab}_{\rm U}-\thetab)
\nonumber
\\[1mm] 
&\!\!=\!\!&
\nub_n^{-1} ({\hat \thetab}_{\rm U}-\thetab)
+
 \lambda(Q\n) \nub_n^{-1}({\hat \thetab}_{\rm C}-{\hat \thetab}_{\rm U} ) 
.
\label{decompP}
\end{eqnarray}
%it is sufficient to show that $\lambda(Q\n) \nub_n^{-1}({\hat \thetab}_{\rm C}-{\hat \thetab}_{\rm U} )$ is $o_{\rm P}(1)$  under~${\rm P}\n_\thetab$.
%, $\thetab\notin\mathcal{M}(\Umb)$. 
For any $\varepsilon >0$, Assumption~A(iii) ensures that 
$$
{\rm P}\n_\thetab[ \lambda(Q\n) \|  \nub_n^{-1} \|  > \varepsilon] \leq {\rm P}\n_\thetab[ \lambda(Q\n)=1]
\to
0, 
$$ 
so that~$\lambda(Q\n) \nub_n^{-1}=o_{\rm P}(1)$ under~${\rm P}\n_\thetab$. Since by assumption, ${\hat \thetab}_{\rm C}-{\hat \thetab}_{\rm U}={\hat \thetab}_{\rm C}-{\thetab}+o_{\rm P}(1)=O_{\rm P}(1)$ under~${\rm P}\n_\thetab$, the result follows from~(\ref{decompP}).
\cqfd
\vspace{4mm}

\noindent {\bf Proof of Theorem~\ref{Vicinityresult}}. 
Writing $
\nub_n^{-1} (\hat \thetab_{\rm U} - \thetab_n) 
=
 \nub_n^{-1} (\hat \thetab_{\rm U} - \thetab) -\nub_n^{-1} (\thetab_n - \thetab) 
%=
% {\bf A}_\thetab {\bf S}_\thetab\n  - \taub_n +o_{\rm P}(1)
%\label{stepeen}
$
and $\nub_n^{-1} (\hat \thetab_{\rm C} - \thetab_n) 
=
\nub_n^{-1} (\hat \thetab_{\rm C} - \thetab) - \nub_n^{-1}(\thetab_n - \thetab) 
%=\Umb{\bf B}_\thetab {\bf S}_\thetab\n  - \taub_n+ o_{\rm P}(1) 
%\nonumber \\
%&=& {\bf A}_\Umb {\bf S}_\thetab\n  - \taub^{(n)}+ o_{\rm P}(1)
%,
%\label{steptwee}
$,
Assumption~B entails that
\begin{equation}
\label{stepdrie}
%{\bf W}_n
%:=
\left(
\begin{array}{c}
\nub_n^{-1} (\hat \thetab_{\rm U} - \thetab_n)\\[-0mm]	
\nub_n^{-1} (\hat \thetab_{\rm C} - \thetab_n)\\[-0mm]	
{\bf D}\n\\[-0mm]
\Lambda\n\\[-2mm] 
\end{array}
\right)
=
\left(
\begin{array}{c}
{\bf A}_\thetab {\bf S}_\thetab\n  - \taub_n\\[-0mm]	
\Umb{\bf B}_\thetab {\bf S}_\thetab\n  - \taub_n\\[-0mm]	
{\bf D}\n\\[-0mm]
\Lambda\n\\[-2mm] 
\end{array}
\right)
+ o_{\rm P}(1) 
\end{equation}
under~${\rm P}\n_\thetab$, $\thetab \in \Thetab_0$. %, hence, from contiguity, also under~${\rm P}\n_{\thetab_n}$. 
Using Assumption~B again, we have
$$
%\left(
%\begin{array}{c}
%\\[-13mm]
%\nub_n^{-1} (\hat \thetab_{\rm U} - \thetab_n)\\[-3mm]	
%\nub_n^{-1} (\hat \thetab_{\rm C} - \thetab_n)\\[-3mm]	
%{\bf D}\n\\[-3mm]
%\Lambda\n\\[-2mm] 
%\end{array}
%\right)
%=
\left(
\begin{array}{c}
{\bf A}_\thetab {\bf S}_\thetab\n  - \taub_n\\[-0mm]	
\Umb{\bf B}_\thetab {\bf S}_\thetab\n  - \taub_n\\[-0mm]	
{\bf D}\n\\[-0mm]
\Lambda\n\\[-2mm] 
\end{array}
\right)
+ o_{\rm P}(1) 
 \stackrel{\mathcal{D}}{\to} 
\mathcal{N}
\left(
\left(
\begin{array}{c}
-\taub\\[-0mm]	
-\taub\\[-0mm]	
{\bf 0}\\[-0mm]
-\frac{1}{2} \taub\pr\Gamb_\thetab \taub\\[-1mm] 
\end{array}
\right)
,
{\bf F}
\right)
$$
under~${\rm P}\n_\thetab$, $\thetab \in \Thetab_0$, with
$$
%\textrm{ and }
{\bf F}
:=
\left(
\begin{array}{cccc} 
{\bf A}_\thetab \Sigb_\thetab {\bf A}_\thetab\pr & {\bf A}_\thetab \Sigb_\thetab   {\bf B}_\thetab\pr \Umb\pr  & {\bf A}_\thetab \Sigb_\thetab {\bf C}_\thetab\pr  & {\bf A}_\thetab \Omegab_\thetab \taub \\[-0mm]
   \Umb{\bf B}_\thetab   \Sigb_\thetab {\bf A}_\thetab \pr &  \Umb{\bf B}_\thetab \Sigb_\thetab {\bf B}_\thetab\pr \Umb\pr &  \Umb{\bf B}_\thetab   \Sigb_\thetab {\bf C}_\thetab \pr & \Umb{\bf B}_\thetab  \Omegab_\thetab \taub \\[-0mm]
    {\bf C}_\thetab \Sigb_\thetab {\bf A}_\thetab\pr &   {\bf C}_\thetab \Sigb_\thetab {\bf B}_\thetab\pr \Umb\pr  & {\bf C}_\thetab \Sigb_\thetab {\bf C}_\thetab\pr  &  {\bf C}_\thetab  \Omegab_\thetab \taub \\[-0mm]
     \taub\pr\Omegab_\thetab {\bf A}_\thetab \pr & \taub\pr \Omegab_\thetab {\bf B}_\thetab\pr \Umb\pr& \taub\pr \Omegab_\thetab {\bf C}_\thetab\pr & \taub\pr \Gamb_\thetab \taub\\[-2mm]
      \end{array} 
\right)
.
$$
%where we used the fact that ${\bf A}_\Umb \Sigb_\thetab {\bf A}_\Umb= {\bf A}_\Umb$ and that ${\bf A}_\Umb \Sigb_\thetab^{1/2} {\bf P}_\Umb={\bf 0}$. 
Thus, the third Le Cam Lemma (jointly with the fact that~(\ref{stepdrie}) also holds under~${\rm P}\n_{\thetab_n}$, from contiguity) directly yields that, under~${\rm P}\n_{\thetab_n}$,
$$
\left(
\begin{array}{c}
\nub_n^{-1} (\hat \thetab_{\rm U} - \thetab_n)\\[-0mm]	
\nub_n^{-1} (\hat \thetab_{\rm C} - \thetab_n)\\[-0mm]	
{\bf D}\n\\[-2mm]
\end{array}
\right)
 \stackrel{\mathcal{D}}{\to} 
\mathcal{N}
\left(
\left(
\begin{array}{c}
(\Ab_\thetab \Omegab_\thetab - {\bf I}_p)\taub\\[-0mm]	
( \Umb  {\bf B}_\thetab\Omegab_\thetab -{\bf I}_p) \taub\\[-0mm]	
{\bf C}_\thetab \Omegab_\thetab \taub \\[-2mm]
\end{array}
\right)
,
\tilde{\bf F}
\right)
,
$$
where~$\tilde{\bf F}$ is obtained from~${\bf F}$ by deleting its last column and last row. 

Therefore, conditional on ${\bf D}\n={\bf D}$, 
$$
\left(
\begin{array}{c}
\nub_n^{-1} (\hat \thetab_{\rm U} - \thetab_n)\\[-0mm]	
\nub_n^{-1} (\hat \thetab_{\rm C} - \thetab_n)\\[-2mm]	
\end{array}
\right)
 \stackrel{\mathcal{D}}{\to} 
\mathcal{N}
\left(
{\bf c}
,
{\bf G}
\right)
$$
under~${\rm P}\n_{\thetab_n}$, where we let
$$
{\bf c}
:=
\left(
\begin{array}{c} 
(\Ab_\thetab \Omegab_\thetab - {\bf I}_p)\taub+ {\bf A}_\thetab \Sigb_\thetab {\bf C}_\thetab\pr  ({\bf C}_\thetab \Sigb_\thetab {\bf C}_\thetab\pr)^- ({\bf D}- {\bf C}_\thetab \Omegab_\thetab \taub) \\[-0mm]
 ( \Umb  {\bf B}_\thetab\Omegab_\thetab -{\bf I}_p) \taub+  \Umb{\bf B}_\thetab   \Sigb_\thetab {\bf C}_\thetab \pr ({\bf C}_\thetab \Sigb_\thetab {\bf C}_\thetab\pr)^- ({\bf D}- {\bf C}_\thetab \Omegab_\thetab \taub) \\[-0mm]
\end{array}
\right)
$$
and
$$
{\bf G}
:=
\left(
\begin{array}{cc}  
{\bf A}_\thetab (\Sigb_\thetab-{\bf L}_\thetab) {\bf A}_\thetab\pr & {\bf A}_\thetab (\Sigb_\thetab-{\bf L}_\thetab)   {\bf B}_\thetab\pr \Umb\pr \\[-0mm]
  \Umb{\bf B}_\thetab (\Sigb_\thetab-{\bf L}_\thetab){\bf A}_\thetab\pr &    \Umb{\bf B}_\thetab (\Sigb_\thetab-{\bf L}_\thetab) {\bf B}_\thetab\pr \Umb\pr  
  \\[-0mm]
  \end{array} 
  \right)
,
$$
with~${\bf L}_\thetab=\Sigb_\thetab {\bf C}_\thetab\pr ( {\bf C}_\thetab \Sigb_\thetab {\bf C}_\thetab\pr)^{-}  {\bf C}_\thetab \Sigb_\thetab$. The result then follows from the fact that $\nub_n^{-1}(\hat \thetab_{\rm PTE} - \thetab_n)=(1- \lambda (\| {\bf D} \|^2)) \nub_n^{-1} (\hat{\thetab}_{\rm U}- \thetab_n)+ \lambda (\| {\bf D} \|^2) \nub_n^{-1} (\hat{\thetab}_{\rm C}- \thetab_n)$ (by using the identities~$\lambda^2 (v)=\lambda (v)$,	$(1-\lambda (v))^2=1-\lambda (v)$, and~$\lambda(v)(1-\lambda(v))=0$). 
\cqfd
\vspace{4mm}

The proof of Proposition~\ref{measqprop} requires the following preliminary result.

\begin{lemma}[\cite{SA06}, pp.~32] 
\label{Salehlem}
Let ${\bf Z}$ be a Gaussian random $p$-vector with mean vector~$\mub$ and covariance matrix~${\bf I}_p$. Then, for any real measurable function $\varphi$, 
$$
(i)
\quad
 {\rm E} [\varphi (\| {\bf Z} \|^2){\bf Z}  ] = {\rm E} [ \varphi ( V)]  \mub
$$
and
$$
(ii)
\quad
{\rm E} [\varphi ( \| {\bf Z} \|^2){\bf Z} {\bf Z}\pr  ] 
=
{\rm E} [ \varphi (V)) ]  {\bf I}_p +   {\rm E} [ \varphi ( W) ]
\mub \mub\pr
,
$$
where $V\sim \chi^2_{p+2}(\| \mub \|^2)$ and $W\sim \chi^2_{p+4}(\| \mub \|^2)$. 
\end{lemma}

\noindent {\bf Proof of Proposition \ref{measqprop}.}  
Since ${\rm E}[{\bf D}]={\bf P}^{\perp}_{\Umb,{\rm eff}} \Gamb_\thetab^{1/2} \taub$ and since ${\bf P}_{\Umb,{\rm eff}}$ is idempotent, we have
\begin{eqnarray}
{\rm E}[{\pmb \mu}_{\rm PTE, eff}^{\rm Vic}] 
&\!=\!& 
{\rm E}[\Gamb_\thetab^{-1/2} {\bf P}^{\perp}_{\Umb,{\rm eff}} ((1- \lambda (\| {\bf D} \|^2)) {\bf D}- \Gamb_\thetab^{1/2} \taub)] 
\nonumber 
\\
&\!=\!&
- \Gamb_\thetab^{-1/2} {\bf P}^{\perp}_{\Umb,{\rm eff}} {\rm E}[\lambda (\| {\bf D} \|^2) {\bf D}]
\label{dd}
.
\end{eqnarray}
Since~${\bf P}_{\Umb,{\rm eff}}^\perp$ is a projection matrix with rank~$p-r$, it decomposes into~${\bf P}_{\Umb,{\rm eff}}^\perp= {\bf O} \Lamb {\bf O}\pr$, where~${\bf O}$ is a $p \times p$ orthogonal matrix and $\Lamb:= {\rm diag}(1, \ldots, 1, 0, \ldots, 0)$ is a diagonal matrix with ${\rm tr}[\Lamb]=p-r$. The random vector~${\bf E}:= {\bf O}\pr {\bf D}$ is then Gaussian with mean vector~$\Lamb {\bf O}\pr \Gamb_\thetab^{1/2} \taub$ and covariance matrix~$\Lamb$. Lemma~\ref{Salehlem}(i) thus entails that 
\begin{equation}
	\label{sal1}
{\rm E}[ \lambda (\| {\bf D} \|^2){\bf D}  ]={\bf O}{\rm E}[ \lambda (\| {\bf E} \|^2){\bf E}  ]=\gamma_2 {\bf P}^{\perp}_{\Umb,{\rm eff}} \Gamb_\thetab^{1/2} \taub
,
\end{equation}
where~$\gamma_2$ is based on a non-central chi-square distribution with~$p-r+2$ degrees of freedom and non-centrality parameter~$(\Lamb {\bf O}\pr \Gamb_\thetab^{1/2} \taub)'\Lamb {\bf O}\pr \Gamb_\thetab^{1/2} \taub=\taub\pr \Gamb_\thetab^{1/2} {\bf P}_{\Umb,{\rm eff}}^\perp \Gamb_\thetab^{1/2} \taub$. 
Plugging this into~(\ref{dd}) provides the result for~${\rm E}[{\pmb \mu}_{\rm PTE, eff}^{\rm Vic}]$.  

We thus turn to~${\rm Var}[{\pmb \mu}_{\rm PTE, eff}^{\rm Vic}]$.  
Since $(1-\lambda(v))^2=1-\lambda(v)$, we have
\begin{eqnarray} 
{\rm Var}[{\pmb \mu}_{\rm PTE, eff}^{\rm Vic}] 
&\!=\!& 
\Gamb_\thetab^{-1/2} {\bf P}^{\perp}_{\Umb,{\rm eff}}  {\rm Var}[(1- \lambda (\| {\bf D} \|^2)) {\bf D}] {\bf P}^{\perp}_{\Umb,{\rm eff}} \Gamb_\thetab^{-1/2} 
\nonumber 
\\
&\!=\!& \Gamb_\thetab^{-1/2} {\bf P}^{\perp}_{\Umb,{\rm eff}} 
\big\{ 
{\rm E}[(1- \lambda (\| {\bf D} \|^2)) {\bf D}{\bf D}\pr] 
\label{stepvars}
  \\
& & 
\hspace{10mm} 
- (1- \gamma_2)^2 {\bf P}^{\perp}_{\Umb,{\rm eff}} \Gamb_\thetab^{1/2} \taub \taub \pr \Gamb_\thetab^{1/2}  {\bf P}^{\perp}_{\Umb,{\rm eff}} 
\big\}
{\bf P}^{\perp}_{\Umb,{\rm eff}} \Gamb_\thetab^{-1/2} 
,
\nonumber
\end{eqnarray}
where we used~(\ref{sal1}). 
Now, by assumption, ${\rm E}[{\bf D}{\bf D}\pr]={\rm Var}[{\bf D}]+{\rm E}[{\bf D}]({\rm E}[{\bf D}])\pr={\bf P}^{\perp}_{\Umb,{\rm eff}}+{\bf P}^{\perp}_{\Umb,{\rm eff}} \Gamb_\thetab^{1/2} \taub\taub\pr \Gamb_\thetab^{1/2} {\bf P}^{\perp}_{\Umb,{\rm eff}}$, and, applying Lemma~\ref{Salehlem}(ii) along the same lines as above, we have that ${\rm E}[\lambda (\| {\bf D} \|^2) {\bf D}{\bf D}\pr]={\bf O} {\rm E}[\lambda (\| {\bf E} \|^2) {\bf E}{\bf E}\pr] {\bf O}'= \gamma_2 {\bf P}^{\perp}_{\Umb,{\rm eff}}+ \gamma_4  {\bf P}^{\perp}_{\Umb,{\rm eff}} \Gamb_\thetab^{1/2} \taub\taub\pr \Gamb_\thetab^{1/2} {\bf P}^{\perp}_{\Umb,{\rm eff}}$. Plugging these expressions into~(\ref{stepvars}) then provides the result.
\cqfd
\vspace{4mm}

\noindent {\bf Proof of Proposition \ref{amseUC}.}  Contiguity implies that~\eqref{optiU} also holds under~${\rm P}\n_{\thetab_n}$, so that
$$
\nub_n^{-1}(\hat \thetab_{\rm U} - \thetab_n) 
= \nub_n^{-1}(\hat \thetab_{\rm U} - \thetab)- \taub 
= \Gamb_\thetab^{-1} \Deltab_\thetab\n- \taub+o_{\rm P}(1)
$$
under~${\rm P}\n_{\thetab_n}$. Since Le Cam's third lemma entails that~$\Deltab_\thetab\n$ is asymptotically normal with mean vector~$\Gamb_\thetab \taub$ and covariance matrix~$\Gamb_\thetab$ under~${\rm P}\n_{\thetab_n}$, it follows that~$\nub_n^{-1}(\hat \thetab_{\rm U} - \thetab_n) $ is asymptotically normal with mean vector~${\bf 0}$ and covariance matrix~$\Gamb_\thetab^{-1}$ under~${\rm P}\n_{\thetab_n}$, which yields~${\rm AMSE}_{\thetab,\taub}(\hat\thetab_{\rm U})=\Gamb_\thetab^{-1}$.
Working along the same lines, we have that, under ${\rm P}\n_{\thetab_n}$, 
\begin{eqnarray}
\nub_n^{-1}(\hat \thetab_{\rm C} - \thetab_n) 
&\!=\!& \nub_n^{-1}(\hat \thetab_{\rm C} - \thetab)- \taub \nonumber \\
&\!=\!& \Umb (\Umb\pr \Gamb_\thetab \Umb)^{-1} \Umb\pr \Deltab_\thetab\n - \taub + o_{\rm P}(1) \nonumber \\
&\!=\!& \Gamb_\thetab^{-1/2} {\bf P}_{\Umb,{\rm eff}}  \Gamb_\thetab^{-1/2}  \Deltab_\thetab\n - \taub + o_{\rm P}(1) \nonumber
.
\end{eqnarray}
It directly follows that $\nub_n^{-1}(\hat \thetab_{\rm C} - \thetab_n)$ is, still under ${\rm P}\n_{\thetab_n}$, asymptotically normal with mean vector $\Gamb_\thetab^{-1/2} {\bf P}_{\Umb,{\rm eff}}  \Gamb_\thetab^{1/2}\taub- \taub= - \Gamb_\thetab^{-1/2}{\bf P}_{\Umb,{\rm eff}}^{\perp} \Gamb_\thetab^{1/2} \taub$ and covariance matrix $\Gamb_\thetab^{-1/2} {\bf P}_{\Umb,{\rm eff}} \Gamb_\thetab^{-1/2}$. The expression for~${\rm AMSE}_{\thetab,\taub}(\hat\thetab_{\rm C})$ given in Proposition~\ref{amseUC} directly follows.
\cqfd

\lhead[\footnotesize\thepage\fancyplain{}\leftmark]{}\rhead[]{\fancyplain{}\rightmark\footnotesize\thepage}

\vspace{4mm}

\noindent 
\vspace{4mm}
\setcounter{section}{1}
\markboth{\hfill{\footnotesize\rm Davy Paindaveine, Jos\'ea  Rasoafaraniaina and Thomas Verdebout} \hfill}
{\hfill {\footnotesize\rm Preliminary test-estimation in LAN models} \hfill}

%
%\bibhang=1.7pc
%\bibsep=2pt
%\fontsize{9}{14pt plus.8pt minus .6pt}\selectfont
%\renewcommand\bibname{\large \bf References}

\bibhang=1.7pc
\bibsep=2pt
\fontsize{9}{14pt plus.8pt minus .6pt}\selectfont
\renewcommand\bibname{\large \bf References}
%\begin{thebibliography}{11}
\expandafter\ifx\csname
natexlab\endcsname\relax\def\natexlab#1{#1}\fi
\expandafter\ifx\csname url\endcsname\relax
  \def\url#1{\texttt{#1}}\fi
\expandafter\ifx\csname urlprefix\endcsname\relax\def\urlprefix{URL}\fi

\vspace{1cm}

%\begin{thebibliography}{11}
%\expandafter\ifx\csname
%natexlab\endcsname\relax\def\natexlab#1{#1}\fi
%\expandafter\ifx\csname url\endcsname\relax
%  \def\url#1{\texttt{#1}}\fi
%\expandafter\ifx\csname urlprefix\endcsname\relax\def\urlprefix{URL
%}\fi
%\bibitem[Antoine and Renault(2012)]{1} Antoine, B. and Renault, E. (2012). Efficient minimum distance
%estimation with multiple rates of convergence. \textit{J. Economet.} \textbf{170}, 350-367.
%\bibitem[Shao and Tu(1995)]{2}
%Shao, J. and Tu, D. (1995). \textit{The Jackknife and Bootstrap}. Springer-Verlag, New York.
%
%%%%%%%%%%%%%%%%%%%%%%%%%%%%%%%%%%%%%%%%%%%%%%%%%%%%%%%%%%%%%%%%%%%%%%%%%%%%%%%%%%%%%%%%%%%%%%%%%%%%%%%%%%%%%%%%%%%%%%%%%%%%%
%\end{thebibliography}

\noindent
ECARES and D\'epartement de Math\'ematique, Universit\'e libre de Bruxelles (ULB)
\vskip 2pt
\noindent
E-mails: dpaindav@ulb.ac.be, rrasoafa@ulb.ac.be, tverdebo@ulb.ac.be

% \vskip .3cm
%\centerline{(Received ???? 20??; accepted ???? 20??)}\par
\end{document}